\documentclass[12pt]{amsart}
\usepackage[latin1]{inputenc}
\usepackage{pdfpages}
\usepackage{amsfonts,amssymb,amsmath,latexsym,stmaryrd,upref}
\usepackage{graphics,amssymb,amsthm,amsmath}
\usepackage[english]{babel} 
\newcommand{\thedate}{\today}

\input xypic \xyoption{all} \CompileMatrices
\newcommand{\pof}{\noindent{\em Proof: }}
\usepackage{color} 
\definecolor{bondiblue}{rgb}{0.0,0.58, 0.71}
\definecolor{cadmiumgreen}{rgb}{0.0, 0.42, 0.24}
\definecolor{coolblack}{rgb}{0.0, 0.18, 0.39}

\newcommand{\ad}{\operatorname{ad}}

\newcommand{\End}{\operatorname{End}}

\newcommand{\Z}{\tiny\operatorname{Z}}

\newtheorem{Thm}{Theorem}[subsection]
\newtheorem{Def}[Thm]{Definition}
\newtheorem{Rem}[Thm]{Remark}
\newtheorem{Obs}[Thm]{Observation}
\newtheorem{Lem}[Thm]{Lemma}
\newtheorem{Cor}[Thm]{Corollary}

\newtheorem{Ex}[Thm]{Example}
\newtheorem{Prop}[Thm]{Proposition} 

\begin{document}

\fontsize{.5cm}{.5cm}\selectfont\sf

\thedate \quad 
Alg-Qdiff-draft.tex

\vskip1cm

\setcounter{tocdepth}{3}

\bigskip \setcounter{section}{0} \bigskip

\title[Algebras of Variable Coefficient Quantized Differential 
Operators]{Algebras of Variable Coefficient Quantized  Differential Operators}

\author{Hans Plesner Jakobsen}

\address{Department of Mathematical Sciences, University of
Copenhagen, Denmark}

\email{jakobsen@math.ku.dk}

\maketitle

\begin{abstract}{In the framework of (vector valued) quantized holomorphic 
functions defined on non-commutative spaces, ``quantized hermitian symmetric 
spaces'',  we analyze  what the algebras of quantized differential operators 
with variable coefficients should be. It is an emediate point that even $0$th 
order operators, given as multiplications by polynomials, have to be specified 
as e.g. left or right multiplication operators since the polynomial algebras 
are replaced by quadratic, non-commutative algebras. In the settings we are 
interested in, there are bilinear pairings which allows us to define 
differential operators as duals of multiplication operators. Indeed, there are 
different choices of pairings which lead to quite different results. We 
consider three different pairings. The pairings are between quantized 
generalized Verma modules and quantized holomorphically induced modules. It is 
a natural demand that the corresponding representations can be expressed by 
(matrix valued) differential operators. We show that a quantum Weyl algebra 
${\mathcal W}eyl_q(n,n)$ introduced by T. Hyashi (\cite{hi}) plays a 
fundamental role. In fact, for one pairing, the algebra of differential 
operators, though inherently depending on a choice of basis, is precisely 
matrices over ${\mathcal W}eyl_q(n,n)$.

We determine explicitly the form of the (quantum) holomorphically induced 
representations and determine, for the different pairings, if they can be 
expressed by differential operators.}
\end{abstract}

\medskip

\medskip

\section{Introduction}

Suppose given 2 quadratic algebras ${\mathcal A}_q^+$ and ${\mathcal A}_q^-$ 
and a non-degenerate bilinear form $(\cdot,\cdot)_X: {\mathcal A}_q^+\times 
{\mathcal A}_q^-\rightarrow {\mathbb C}$. We will say that $(\cdot,\cdot)_X$ 
gives a pairing between the two algebras; the label $X$ reflects 
that there will be several different pairings.  A pairing displays the 
algebras as duals of each other, or more correctly, sets up identifications 
between one algebra and the dual of the other. While the dual space is unique, 
there may be considerable interest and usefulness in exhibiting the duals 
concretely in such manners. Even when the two algebras are given in advance, 
one may examine different pairings between them to optimize certain properties.

\medskip

An algebra $A$ is of course a left and a right module over itself. 
Thus, one has left multiplication operators ${}_aM$ and right multiplication 
operators $M_b$ for all $a,b\in A$. Clearly, the algebra generated by the 
operators ${}_aM$ is isomorphic to $A$ while  the algebra determined by the 
operators $M_b$ is  isomorphic to $A^{o}$; the opposite algebra.

\medskip

In the case with two algebras in a pairing such as the previously mentioned 
algebras ${\mathcal A}_q^\pm$, we can define left and right ``constant 
coefficient'' differential operators ${}_c\partial_X,{}_X\partial_c$, on, say, 
${\mathcal A}_q^+$, as the operators obtained as the duals, via the pairing, of 
left and right multiplication operators on ${\mathcal A}_q^-$. These operators, 
which depend upon the pairing $X$, may then be put together with the left and 
right multiplication operators to form an algebra ${\mathcal D}^+_{X, Full}$ of 
differential operators. Instances of such (left) operators are the Kashiwara 
derivations (\cite{K1}). In specific examples it is  interesting to determine 
if ${\mathcal D}^+_{X, Full}$, for a specific index $X$, is generated by fewer 
operators, eg by left differential operators and left multiplication operators; 
${\mathcal D}^+_{X, Full}=Alg_{\mathbb C}\left\{{}_dM, {}_c\partial_X\; ; c\in 
{\mathcal A}_q^-, d\in {\mathcal A}_q^+\right\}$?

\medskip

The algebras ${\mathcal A}_q^\pm$ we  consider are quadratic algebras that are 
specific subalgebras of ${\mathcal U}_q(su(n,n)^{\mathbb C})$. As quadratic 
algebras they are actually isomorphic. They are the algebras of the quantized 
generalized unit disk. Furthermore, they are modules for ${\mathcal 
U}_q({\mathfrak k}^{\mathbb C})$, where ${\mathfrak k}$ is a maximal compact 
subalgebra of $su(n,n)$, but as such, they are non-isomorphic; and not 
necessarily dual modules, either. 

\medskip

There is a further structure we need to 
include in our discussions: To each finite-dimensional 
${\mathcal U}_q({\mathfrak k}^{\mathbb C})$ module $V_\Lambda$ there is a 
quantized generalized Verma module over ${\mathcal U}_q(su(n,n)^{\mathbb C})$. 
As a ${\mathcal U}_q({\mathfrak k}^{\mathbb C})$ module it is given as 
${\mathcal M}(V_\Lambda)= {\mathcal A}_q^-\otimes V_\Lambda$.

This extra structure leads to natural demands on the bilinear pairing. In 
this connection it is not profitable to consider a pairing between a Verma 
module and its ``opposite'' (interchanging positive and negative roots). The 
right notion of a dual of a generalized Verma module, in our context, is a 
holomorphically induced module.

\medskip

It is natural to demand that an algebra of differential operators is rich 
enough that the operators in the holomorphically induced modules belong to it. 
To have any hope of that,  one will of course need to include in the algebra 
the homomorphisms $hom_{\mathbb C}(V_\Lambda,V_\Lambda)$, or, rather, the duals 
thereof.

\medskip

About pairings: In the classical situation, the Killing form on a real 
semi-simple Lie algebra ${\mathfrak g}$, extended to ${\mathcal U}({\mathfrak 
g}^{\mathbb C})\times {\mathcal U}({\mathfrak g}^{\mathbb C})$, gives the 
wanted pairing. M Rosso  constructed the quantum analogue of this. After that, 
G. Lusztig (\cite{luz2}) and M. Kashiwara (\cite{K1}) made valuable extensions 
and simplifications, and Kashiwara defined some derivations as duals of left or 
right multiplications. The history of this subject is very rich and we hope 
that we are not being too unfair in this sketchy summary. One should definitely 
also consult \cite{jo-le} and \cite{tani}. We follow here the book by J.C. 
Jantzen (\cite{jan}), not only for notation, but actually  to the extent of 
copying directly several of his constructions and results. 
\medskip

We will study three bilinear forms, indexed by $X=J, K, L$. $(\cdot,\cdot)_J$ 
is the form considered by Jantzen, though he actually studies an additional 
form towards the end of his book. They reflect the three standard ways of 
quantizing integers:

\begin{eqnarray*}[[a]]_q&=&1+q^2+\cdots + q^{2a-2}\ (J),\\
\{\{a\}\}_q&=&1+q^{-2}+\cdots + q^{-(2a-2)}\ (K),  \textrm{ and
}\\
\left[a\right]_q&=&q^{-a+1}+\cdots+q^{a-1}\ (L).\end{eqnarray*}

\medskip

Using fixed PBW bases, we define an auxiliary algebra ${\mathcal W}eyl_q(n,n)$ 
- a quantization of the classical  Weyl algebra in $n^2$ variables - as the 
algebra generated by $n^2$ commuting variables $D_{i,j}$ and $n^2$ commuting 
variables $M_{i,j}$, and where also $D_{i,j}$ commutes with $M_{s,t}$ if 
$(i,j)\neq (s,t)$ so that the only non-trivial relations are at fixed nodes. 
Here the relations are\begin{eqnarray*}\nonumber
D_{i,j}M_{i,j}-qM_{i,j}D_{i,j}=H_{i,j}^{-1},&&
D_{i,j}M_{i,j}-q^{-1}M_{i,j}D_{i,j}=H_{i,j},
\\\nonumber H_{i,j}D_{i,j}=q^{ -1}
D_{i,j}H_{i,j},\textrm{ and }  &&H_{i,j}M_{i,j}
=qM_{i,j}H_{i,j} .\end{eqnarray*}

This is a very interesting algebra which was introduced by T. Hyashi 
(\cite{hi}). There has recently been renewed interest in it, see \cite{kilu}.

\medskip

It turns out that there is a big difference between the three cases, where 
especially the case (J) leads to unpleasant  results. The simplest case, on the 
other hand, is the case (L) where 
\begin{eqnarray*}Alg_{\mathbb C}\left\{{}_dM, {}_c\partial_L\; ; c\in {\mathcal 
A}_q^-, d\in {\mathcal A}_q^+\right\}&=&Alg_{\mathbb C}\left\{M_d, 
{}_L\partial_c\; ; c\in {\mathcal A}_q^-, d\in {\mathcal 
A}_q^+\right\}\\&=&{\mathcal W}eyl_q(n,n).\end{eqnarray*}

For the case (K) there is a big subalgebra ${\mathcal K}{\mathcal W}eyl_q(n,n)$ 
of ${\mathcal W}eyl_q(n,n)$ with many pleasing properties such that 

\begin{eqnarray*}&&Alg_{\mathbb C}\left\{{}_dM ,M_f ,{}_c\partial_K\; ; c\in 
{\mathcal A}_q^-, d, f\in {\mathcal A}_q^+\right\}=\\&&Alg_{\mathbb 
C}\left\{M_d, {}_fM,  {}_K\partial_c\; ; c\in {\mathcal A}_q^-, d,f \in 
{\mathcal A}_q^+\right\}\ =\ {\mathcal K}{\mathcal W}eyl_q(n,n).\end{eqnarray*}

In both of the cases (K), (L), these algebras, augmented by (constant value) 
matrices,  contain the generators of the holomorphically induced 
representations. To prove such a statement it suffices to determine the action 
of $E_\beta$, where  $\beta$ is the unique non-compact root and prove the 
statement in this special case. The mentioned action is given in 
Corollary~\ref{actJ} to Theorem~\ref{main1}; our first main result.

\medskip

There is one related important study, namely that by L. Vaksman and his group 
\cite{vaks}. In it, they extend substantially the quantized exterior derivative 
introduced in \cite{chari} and which already leads to derivatives. Vaksman et 
al.  discovered a fundamental symmetry. Their method uses induction from the 
trivial ${\mathcal U}_q({\mathfrak k}^{\mathbb C})$ module. Here, the 
classically holomorphically induced module is annihilated by first order 
differential operators. In the quantized situation, there is a natural pairing 
and a natural algebra structure  obtained from the tensor product. Proceeding 
like this, they have frozen the algebras at a specific weight which means that 
they do obtain interesting results, but not the general picture that we obtain. 
The extra symmetry is related to the fact that our algebras are bi-modules. 
Furthermore, their pairing is degenerate, though ``mildly''.

\smallskip

Quantized differential operators were also studied in \cite{j-qdiff}, 
\cite{j-qdirac}, \cite{genkai}.

\medskip

In \S2 we introduce the quantized Hermitian symmetric spaces (the case of 
$su(n,n)$) via the Lusztig operators. In \S3 we study these spaces as 
${\mathcal U}_q({\mathfrak k}^\mathbb C)$ modules and for this purpose, 
${\mathcal W}eyl_q(n,n)$ is introduced. \S4 contain many direct quotes from 
Jantzen's book. The bilinear pairing is introduced, the duals to left and right 
multiplication operators are determined, and the left action of $E_\alpha$ 
($\alpha$ a simple root)  in ${\mathcal U}_q^-$ is given. Duality considerations 
are continued in \S5 where the various pairings we wish to study, are 
introduced. We also discuss various change-of-basis maps. One such is needed 
because Jantzen's form has a singularity at $q=1$. In \S6 we introduce the 
generalized Verma modules and the quantized holomorphically induced modules, and 
pairings between them. 

Then, in \S7 we obtain the dual of the action of $E_\beta$, $\beta$ the unique 
non-compact simple root. This is given in Theorem~\ref{main1}. It should 
be observed how simple the result actually is and that it is given, essentially, 
by left and right multiplication operators. We also obtain the limit of the 
operator as $q\rightarrow1$ and make sure that it agrees with the known 
``classical'' operator.

 Finally, in \S8 we obtain the other main results about the algebras of 
polynomial coefficient differential operators; Theorem~\ref{L-thm} and 
Theorem~\ref{K-thm}. As a bonus we obtain that ${\mathcal W}eyl_q(n,n)$, though 
manifestly defined via a PBW basis, actually is intrinsic.
 
The main technical part involves computing explicitly the multiplication 
operators ${}_cM,M_d$ which can be done using the defining quadratic relations. 
If we let $D_{i,j}$ and $M_{i,j}$ denote the generators of ${\mathcal 
W}eyl_q(n,n)$, $i,j\in\{1,2,\dots,n\}$, and similarly let $H_{i,j}^{\pm1}$ 
denote the operators $H^{\pm1}$ at the node $i,j$ then in all cases $(J), (K), 
(L)$ the following holds: At each node we obtain from the left and right 
operators along with their duals, operators $D_{i,j}\psi_{i,j}$ and 
$M_{i,j}\phi_{i,j}$ for some elements $\psi_{i,j}$ and $\phi_{i,j}$ which are 
Laurent monomials in the elements $H_{1,1}, H_{2,1},\dots, H_{n-1,n},H_{n,n}$. 
The appearance of these factors is just one of the interesting consequences of 
working with a quantized Weyl algebra.

In the cases (L) and (K) we get sufficiently many such elements to find some 
simple generators. In the case (J), however, the generators remain complicated.

\bigskip

\section{Quantized Hermitean Symmetric spaces.}

\subsection{Basic definitions.}
	  
We consider ${\mathfrak g}=su(n,n)$; ${\mathfrak g}^{\mathbb C}$ is a simple 
complex Lie algebra 
of type $A_{2n-1}$. We choose below a set of simple roots $\Pi$
in the root space $\Psi$. The Weyl group is equal to
$S_{2n}$. We denote the generators of the Weyl group by $s_\gamma, 
\gamma\in\Pi$ and denote by
$E_\gamma,F_\gamma,K_\gamma^{\pm1}$ for $\gamma\in\Pi$ the generators of 
${\mathcal U}_q({\mathfrak g}^{\mathbb C})$ is standard notation. The weight 
lattice is 
denoted by ${\mathfrak L}$ and we further extend the notation $K_\xi$
to hold for any weight $\xi \in {\mathfrak L}$ in the usual fashion.

\smallskip

The roots $\Psi$ may be represented in ${\mathbb
R}^{2n}$ by the set \begin{equation}\Psi=\{\pm e_i\mp e_j\mid
i,j=1,2,\dots, 2n\textrm{ and }i\neq j\},\end{equation} where
$\{e_1,e_2,\dots, e_{2n}\}$ is the standard basis of
${\mathbb R}^{2n}$. We then have
\begin{equation}\Pi=\{e_i-e_{i+1}\mid i= 1,2,\dots,
2n-1\}.\end{equation} Throughout, we let
$\beta=e_n-e_{n+1}$ denote the unique {\bf non-compact} simple
root. The roots $\nu_i=e_{n-i}-e_{n-i+1}, i=1,\dots, n-1$ and
the roots $\mu_j=e_{n+j}-e_{n+j+1}, j=n+1,\dots, 2n$ are the {\bf compact}
simple roots of type $A_{n-1}$; $\Pi=\{\mu_1,\dots,\mu_{n-1}\}\cup\{\beta\}\cup
\{\nu_1,\dots,\nu_{n-1}\}$.  We also set $\Pi_c=\{\mu_1,\dots,\mu_{n-1}\}\cup
\{\nu_1,\dots,\nu_{n-1}\}$; the compact simple roots, and set 
$\Pi_L=\{\mu_1,\dots,\mu_{n-1}\}$, $\Pi_R=
\{\nu_1,\dots,\nu_{n-1}\}$. Let ${\mathfrak 
k}_L^{\mathbb C}$ and  ${\mathfrak k}_R^{\mathbb C}$ denote the subalgebras 
defined by the simple roots  $\{\mu_1,\dots,\mu_{n-1}\}$ and 
$\{\nu_1,\dots,\nu_{n-1}\}$, respectively. Finally, for $k=1,\dots, n-1$, let 
$\delta^\mu_k$ and 
$\delta^\nu_k$, respectively, denote the fundamental dominant weights for the 
roots $\{\mu_1,\dots,\mu_{n-1}\}$ and $\{\nu_1,\dots,\mu_{n-1}\}$, respectively.

\bigskip

A maximal compact subalgebra ${\mathfrak k}^{\mathbb
C}$ of ${\mathfrak g}^{\mathbb
C}$ is
\begin{equation}
{\mathfrak k}^{\mathbb
C}=su(n)^{\mathbb C}\oplus {\mathbb C}\oplus
su(n)^{\mathbb C}={\mathfrak k}_L^{\mathbb C}\oplus
{\zeta} \oplus {\mathfrak k}_R^{\mathbb C},
\end{equation}where $\zeta$ is the center of ${\mathfrak k}^{\mathbb
C}$ and is generated by an element $h_\beta$ of the compact Cartan subalgebra. 
We have furthermore on the classical level

\begin{equation}{\mathfrak
g}^{\mathbb C}={\mathfrak
k}^{\mathbb C}\oplus {\mathfrak p}={\mathfrak p}^-\oplus
{\mathfrak k}^{\mathbb C}\oplus {\mathfrak
p}^+,
\end{equation} where ${\mathfrak
p}^{\pm}$ are abelian ${\mathfrak k}^{\mathbb C}$ 
modules, and \begin{equation}{\mathcal U}({\mathfrak
g}^{\mathbb C})={\mathcal P}({\mathfrak p}^-)\cdot
{\mathcal U}({\mathfrak k}^{\mathbb C})\cdot {\mathcal
P}({\mathfrak p}^+).\label{class}\end{equation} We have that ${\mathcal 
U}({\mathfrak p}^\pm)={\mathcal P}({\mathfrak p}^\pm)$ are polynomial algebras.

On the quantized level we have
\begin{eqnarray}  {\mathcal
U}_q({\mathfrak k}^{\mathbb C})={\mathcal
U}_q({\mathfrak k}_L^{\mathbb C})\cdot {\mathbb
C}[K_\beta^{\pm1}] \cdot {\mathcal U}_q({\mathfrak
k}_R^{\mathbb C}). \end{eqnarray} 
There is an analogue of (\ref{class}) for ${\mathcal
U}_q({\mathfrak g}^{\mathbb C})$,
\begin{equation} {\mathcal
U}_q({\mathfrak g}^{\mathbb C})={\mathcal A}_q^-\cdot
{\mathcal U}_q({\mathfrak k}^{\mathbb C})\cdot {\mathcal
A_q}^+.\end{equation} Here, ${\mathcal A}_q^\pm$ are quadratic
algebras which are furthermore ${\mathcal U}_q({\mathfrak
k}^{\mathbb C})$. We will describe these later.

We let ${\mathcal U}^0_q({\mathfrak k}_L^{\mathbb C}\oplus {\mathfrak 
k}_R^{\mathbb C})$  denote the Laurent polynomials  generated by the elements 
$K_\alpha$ for $\alpha\in\Pi_c$, and let ${\mathcal U}^0_q({\mathfrak 
k}^{\mathbb C})$ denote the analogue where also $\alpha=\beta$ is allowed.

\bigskip

For use in the construction of the algebras ${\mathcal A}_q^\pm$ we now 
consider some elements in the Weyl group. 

Let $I=(1,2,\dots,n,n+1,\dots,2n)$. Consider the
following elements in $S_{2n}$: \begin{eqnarray}
\omega_0(I)&=&(2n,2n-1,\dots,n+1,n,\dots,2,1)\\
\omega^+(I)&=&(n+1,\dots,2n-1,2n,1,2,\dots,n-1, n)\\\label{w-mu}
\omega^L(I)&=&(n,n-1,\dots,2,1,n+1,n+2,\dots,
2n-1,2n)\\
\omega^R(I)&=&(1,2,\dots,n-1,n,2n,2n-1,\dots,n+2
,n+1). \end{eqnarray}

We have that $\omega_0$ is the longest element, and
\begin{equation}\omega_0=\omega_0^{L}
\omega_0^{R}\omega_0^+=\omega_0^+\omega_0^{L}
\omega_0^{R}=
\omega_0^{L}\omega_0^+\omega_0^{L},\end{equation}
and many more similar identities. We shall later use
\begin{equation}\omega_0=\omega_0^{L}
\omega_0^{R}\omega_0^+\end{equation} since this puts
the  elements from $\Pi_c$ to the right according to
the construction in \cite[p. 163 -168]{jan}. However, it is
convenient first to consider
$\omega_0=\omega_0^+\omega_0^{L} \omega_0^{R}$:

We shall in all cases use the following choice for a
reduced expression for $\omega_0^+$:
\begin{equation}\omega_0^+=s_\beta s_{\mu_1}\cdots
s_{\mu_{n-1}}s_{\nu_1}s_\beta s_{\mu_1}\cdots
s_{\mu_{n-2}}s_{\nu_2}s_{\nu_1}s_\beta \cdots s_\beta
s_{\mu_1} s_{\nu_{n-1}}\cdots s_{\nu_1}s_\beta.\end{equation}

\begin{Rem}If we look at $n\times(n+r)$ we replace $\omega_0^+$ by 
$\omega_0^+\omega_E^+$ where \begin{equation}\omega_E^+=s_{\nu_{n}}\cdots 
s_{\nu_1}s_{\nu_{n+1}}\cdots s_{\nu_2}\dots s_{\nu_{n+r-1}}\cdots 
s_{\nu_r}.\end{equation}\end{Rem}

\medskip

\subsection{The Lusztig operators}

G. Lusztig (\cite{l}) has given a construction of braid
operators $T_\alpha$ for each simple root $\alpha$, and
has extended these to operators $T_\omega$ for each
element $\omega$ of the Weyl group. Indeed, there are
two choices of such operators, the other usually denoted
by $T'_\omega$, and there is always a choice between $q$ and
$q^{-1}$. In the book (\cite{jan}), J. C. Jantzen describes, among many other 
things, these operators. We will throughout use the notation and choices from 
this book.
  (\cite{jan}).

\medskip

\begin{Lem}\label{156}[\cite{l}, \cite{jan}[p. 156]] In the simply
laced case with neighboring simple roots $\alpha,\gamma$:
\begin{eqnarray} \label{9}T_\alpha (F_\gamma)&=&F_\gamma
F_\alpha-qF_\alpha F_\gamma\\ T_\alpha(E_\gamma)&=&E_\alpha
E_\gamma-q^{-1} E_\gamma E_\alpha=\ad(E_\alpha)(E_\gamma)\\
T^{-1}_\alpha (F_\gamma)&=&F_\alpha F_\gamma - qF_\gamma F_\alpha\\
T^{-1}_\alpha (E_\gamma)&=&E_\gamma E_\alpha - q^{-1}E_\alpha
E_\gamma. \end{eqnarray} \end{Lem}
	
	\medskip

\begin{Lem}\label{157}[\cite{l}, \cite{jan}[p. 156]] If $\alpha,\gamma$ are 
adjacent, 
then\begin{eqnarray}
T_\gamma(E_\alpha)&=&T^{-1}_\alpha(E_\gamma)\\
T_\gamma(F_\alpha)&=&T^{-1}_\alpha(F_\gamma)\\ T_\alpha T_\gamma
T_\alpha&=& T_\gamma\label{lu} T_\alpha T_\gamma \\T_\alpha
T_\gamma(E_\alpha)&=&E_\gamma\\T_\alpha
T_\gamma(F_\alpha)&=&F_\gamma. \end{eqnarray}
If $\alpha,\gamma$ are not adjacent, then 
\begin{equation}
 T_\alpha T_\gamma= T_\gamma T_\alpha
\end{equation}
\end{Lem}

\smallskip

Lusztig has further shown that if $s_\alpha s_\beta$ has
order $m$ then \begin{equation}\underbrace{T_\alpha T_\gamma \cdots}_m= 
\underbrace{T_\gamma
T_\alpha\cdots}_m.\end{equation}

\smallskip

There is an important construction of a PBW type
basis in ${\mathcal U}^{\pm}$ for any reduced
decomposition of $\omega_0=s_{\alpha_1}s_{\alpha_2}\dots
s_{\alpha_{i}}\dots s_{n(2n-1)}$. It is also due to Lusztig  and is given as
follows: Set, $\omega_i=s_{\alpha_1}s_{\alpha_2}\dots s_{\alpha_{i}}, 
i=1,\dots,2n-1$.
Then, \begin{eqnarray}\label{L1}
\gamma_i&=&\omega_{i-1}(\alpha_i)\\\label{L2}
X_{\gamma_i}&=&T_{\omega_{i-1}}(E_{\alpha_i})=T_{\alpha_1}T_{\alpha_2}\dots
T_{\alpha_{i-1}}(E_{\alpha_i})\in U^+_{\gamma_i}\\
Y_{\gamma_i}&=&T_{\omega_{i-1}\label{L3}}(F_{\alpha_i})=T_{\alpha_1}T_{\alpha_2}
\dots
T_{\alpha_{i-1}}(F_{\alpha_i}) \in
U^-_{\gamma_i}\end{eqnarray}

	  \medskip
	  
We now introduce some intermediary bases. Set
\begin{eqnarray}Q_{ij}&=&T^{-1}_{\mu_{i-1}}T^{-1}_{\mu_{i-2}}
\dots
T^{-1}_{\mu_{0}}T^{-1}_{\nu_{j-1}}T^{-1}_{\nu_{j-2}}\dots
T^{-1}_{\nu_{0}}(F_\beta)\\
P_{ij}&=&T^{-1}_{\mu_{i-1}}T^{-1}_{\mu_{i-2}}\dots
T^{-1}_{\mu_{0}}T^{-1}_{\nu_{j-1}}T^{-1}_{\nu_{j-2}}\dots
T^{-1}_{\nu_{0}}(E_\beta). \end{eqnarray}
	
Set ${\bf a}=(a_{11},\dots, a_{1,n}\dots,
a_{nn})\in{\mathbb N}_0^{n^2}$, and set
\begin{eqnarray}P^{\mathbf a}&=&P_{11}^{a_{11}}\cdots
P_{nn}^{a_{nn}}\\ Q^{\mathbf b}&=&Q_{11}^{b_{11}}\cdots
Q_{nn}^{b_{nn}}. \end{eqnarray}
	
Let $\{u^{{\mathfrak k}^+}_i\mid i\in I\}$ denote a PBW
type basis of ${\mathcal U}^+({\mathfrak k})$ and let
$\{u^{{\mathfrak k}^-}_j\mid j\in J\}$ denote a PBW type
basis of ${\mathcal U}^-({\mathfrak k})$.
	
\begin{Prop}There is a basis \begin{equation}\{u^{{\mathfrak
k}^+}_i\cdot P^{\mathbf a}\mid i\in I; {\bf
a}\in{\mathbb N}_0^{n^2}\}\end{equation}of ${\mathcal U}^+$, and a
basis \begin{equation}\{u^{{\mathfrak k}^-}_j\cdot Q^{\mathbf b}\mid
j\in J; {\bf b}\in{\mathbb N}_0^{n^2}\}\end{equation}of ${\mathcal
U}^-$. \end{Prop}

\proof This follows from Lusztig (see Jantzen \S8.24) by
using the following extra observations which are easily
deduced from, in particular, Lemma~\ref{157}.
	
1) $T_\beta T_{\mu_1}(E_{\mu_2})=T_\beta
T^{-1}_{\mu_2}(E_{\mu_1})=T^{-1}_{\mu_2}T_\beta(E_{\mu_1})=T^{-1}_{\mu_2}T^{-1}_
{\mu_1}(E_\beta)$.

\smallskip

2) $T_\beta
T_{\mu_1}T_{\mu_2}T_{\nu_1}(E_\beta) =T_\beta
T_{\mu_1}T_{\nu_1}(E_\beta)=T_\beta
T_{\mu_1}T_\beta^{-1}(E_{\nu_1})
=$

$T_{\mu_1}^{-1}T_{\beta}T_{\mu_1}(E_{\nu_1})=
T_{\mu_1}^{-1}
T_{\beta}(E_{\nu_1})=T_{\mu_1}^{-1}T_{\nu_1}^{-1}(E_\beta).
$\qed

\medskip

While these bases have many good properties, it is more natural to have 
${\mathfrak k}^+$ to the
right. To this end we employ $\omega_0=\omega^L\omega^R \omega^+$:
	
	\smallskip

\begin{Def}Set, for all $i,j\in\{1,2,\dots,n\}$,
\begin{equation}Z_{n+1-i,n+1-j}=T_{\omega^\mu\omega^\nu}(P_{ij})\textrm{
and }W_{n+1-i,n+1-j}=T_{\omega^\mu\omega^\nu}(Q_{ij}).\end{equation}
\end{Def}

Let the roots $\gamma_{ij}$ be defined by 
\begin{equation}\forall i,j\in\{1,2,\dots,n\}: Z_{ij}=Z_{\gamma_{ij}}\textrm{ 
and
}W_{ij}=W_{-\gamma_{ij}}.\end{equation}

\medskip
\begin{Lem}\begin{eqnarray}Z_{i,j}&=&T_{\nu_{j-1}}
T_{\nu_{j-2}}\dots T_{\nu_{0}}\cdot T_{\mu_{i-1}}
T_{\mu_{i-2}}\dots T_{\mu_{0}}(E_\beta)\\
W_{i,j}&=&T_{\nu_{j-1}} T_{\nu_{j-2}}\dots
T_{\nu_{0}}\cdot T_{\mu_{i-1}} T_{\mu_{i-2}}\dots
T_{\mu_{0}}(F_\beta).\label{26}\end{eqnarray} In particular,
\begin{equation} Z_{i,j}=ad(E_{\nu_{j-1}})
ad(E_{\nu_{j-2}})\dots ad(E_{\nu_{0}})\cdot ad(E_{\mu_{i-1}})
ad(E_{\mu_{i-2}})\dots ad(E_{\mu_{0}})(E_\beta).
\end{equation}
\end{Lem}
	  
	  \smallskip
	  
\proof 
Since all steps are similar, it suffices to prove that
\begin{equation}T_{\omega^\mu}(T^{-1}_{\mu_{i}}T^{-1}_{\mu_{i-1}}\dots
T^{-1}_{\mu_{0}}(E_\beta))=T_{\mu_{n-i}}T_{\mu_{n-i-1}}\dots
T_{\mu_{0}}(E_\beta).\end{equation} Using
\begin{equation}(\omega^L)=s_{\mu_{1}}s_{\mu_{2}}
s_{\mu_{1}}s_{\mu_{3}}s_{\mu_{2}}s_{\mu_{1}}\dots
s_{\mu_{n-1}}s_{\mu_{n-2}}\dots s_{\mu_{1}},\end{equation}this
follows from the above formulas, especially (\ref{lu}) in
the form \begin{equation}T_{\mu_{k-1}}
T_{\mu_{k}}^{-1}T_{\mu_{k-1}}^{-1}=
T_{\mu_{k}}^{-1}T_{\mu_{k-1}}^{-1} T_{\mu_{k}}\end{equation} and
$T_{\mu_{k}}^{\pm1}(X_\beta)=0$ for $k\geq2$ and $X=E$
or $X=F$. \qed

\medskip

Observe that we now have an ordering 
\begin{eqnarray}\label{long-order}W_{1,1}, W_{2,1},\dots,W_{n,1}W_{1,2},\dots, 
W_{n,2},\dots,W_{1,n},\dots, W_{n,n},\\\dots W_{1,n+1},\dots, W_{n,n+1},\dots, 
W_{1,n+r},\dots,W_{n,n+r},\nonumber\end{eqnarray}
and a similar ordering of the elements $Z_{i,j}$.
	  
	  \bigskip
	  
	  Now define 
	  
\begin{eqnarray}Z^{\mathbf
a}&:=&Z_{11}^{a_{11}}Z_{21}^{a_{21}}\cdots
Z_{n1}^{a_{n1}}Z_{12}^{a_{12}}\cdots
Z_{nn}^{a_{nn}},\label{20} \\W^{\mathbf
a}&:=&W_{11}^{a_{11}}W_{21}^{a_{21}}\cdots
W_{n1}^{a_{n1}}W_{12}^{a_{12}}\cdots\label{21}
W_{nn}^{a_{nn}}.\end{eqnarray}

\medskip
	 
\begin{Prop}\label{2.5}There is a basis \begin{equation}\{ Z^{\mathbf
a}\cdot u^{{\mathfrak k}^+}_i\mid i\in I; {\bf
a}\in{\mathbb N}_0^{n^2}\}\end{equation}of ${\mathcal U}^+$, and a
basis \begin{equation}\{ W^{\mathbf b}\cdot u^{{\mathfrak k}^-}_j\mid
j\in J; {\bf b}\in{\mathbb N}_0^{n^2}\}\end{equation}of ${\mathcal
U}^-$. \end{Prop}

\smallskip
	 
These are the bases we will use.

	  \bigskip
	    
\begin{Ex}$su(3,3)$:\begin{equation}\omega_0=(\beta\mu_1\mu_2
\nu_1\beta\mu_1\nu_2\nu_1\beta)(\mu_1\mu_2\mu_1)
(\nu_2\nu_1\nu_2).\end{equation}

Then,{\small
\begin{equation}\begin{array}{cccccccccccccc}\gamma_{11}&=&\beta,&
\gamma_{12}&=&s_\beta
s_{\mu_1}s_{\mu_2}({\nu_1}),&\gamma_{13}&=&
s_{\beta}s_{\mu_1}s_{\mu_2} s_{\nu_1}s_\beta
s_{\mu_1}({\nu_2}),\\
\gamma_{{21}}&=&s_\beta({\mu_1}),&
\gamma_{{22}}&=&s_{\beta}s_{\mu_1}s_{\mu_2}
s_{\nu_1}(\beta),&\gamma_{{23}}&
=&s_{\beta}s_{\mu_1}s_{\mu_2} s_{\nu_1}s_\beta
s_{\mu_1}s_{\nu_2}({\nu_1}),\\
\gamma_{31}&=&s_{\beta}s_{\mu_1}({\mu_2}),&\gamma_{32}&
=&s_{\beta}s_{\mu_1}s_{\mu_2}
s_{\nu_1}s_\beta({\mu_1}),&\gamma_{{33}}&
=&s_{\beta}s_{\mu_1}s_{\mu_2} s_{\nu_1}s_\beta
s_{\mu_1}s_{\nu_2}s_{\nu_1}(\beta)\ . \end{array}\end{equation}}
	  
	\medskip
	
	Furthermore,
	{\small
\begin{equation}\begin{array}{ccc} Z_{11}=E_\beta,
&Z_{12}=E_{\nu_1}Z_{11}-q^{-1}Z_{11}E_{\nu_1},&Z_{13}=E_{\nu_2}Z_{12}-q^{-1}Z_{
12
}E_{\nu_2},\\
Z_{21}=E_{\mu_1}Z_{11}-q^{-1}Z_{11}E_{\mu_1},&Z_{22}=E_{\nu_1}Z_{21}-q^{-1}Z_{21
}
E_{\nu_1},&Z_{23}=E_{\nu_2}Z_{22}-q^{-1}Z_{22}E_{\nu_2},\\
Z_{31}=E_{\mu_2}Z_{21}-q^{-1}Z_{21}E_{\mu_2},&Z_{32}=E_{\nu_1}Z_{31}-q^{-1}Z_{31
}
E_{\nu_1},&Z_{33}=E_{\nu_2}Z_{32}-q^{-1}Z_{32}E_{\nu_2}\ .
\end{array}\end{equation}}

It follows easily from the quantized Serre relations
that these elements generated the usual quantized
$3\times3$ matrix algebra with ``$q^{-1}$ relations'' (see the Definition 
below).\end{Ex}

\medskip

\subsection{The quadratic algebras ${\mathcal A}_q^+, {\mathcal A}_q^-$.
}	  

\begin{Def}We let ${\mathcal A}_q^+$ and ${\mathcal
A}_q^-$ denote the algebras generated by the elements 
$Z_{ij},i,j\in\{1,2,\dots,n\}$ and $W_{ij},i,j\in\{1,2,\dots,n\}$, 
respectively.\end{Def}
	
	\medskip

\begin{Prop}\label{qrel}$\forall i,j,k,s,t\in\{1,2,\dots,n\}$:
\begin{eqnarray}\label{a}Z_{ij}Z_{ik} &=&
q^{-1}Z_{ik}Z_{ij} \textrm{ if }j < k;\\
\label{b}Z_{ij}Z_{kj} &=& q^{-1}Z_{kj}Z_{ij}\textrm{ if
}i< k;\\\label{c} Z_{ij}Z_{st} &=& Z_{st}Z_{ij} \textrm{
if }i < s\textrm{ and }t < j;\\\label{cross}
Z_{ij}Z_{st} &=& Z_{st}Z_{ij}-(q-q^{-1})Z_{it}Z_{sj} =
\textrm{ if }i < s\textrm{ and } j < t. \end{eqnarray}
There are entirely identical relations for the elements $W_{ij}$.
\end{Prop}	
	
\proof 	These relations follow from Lemma~\ref{156}, Lemma~\ref{157}, and the 
quantum Serre relations. \qed

	\medskip

\begin{Rem} We see that one  gets the relations studied in e.g. (\cite{jz}) 
except for the replacement 
$q\rightarrow q^{-2}$. In  (\cite{jz}) the methods were not related to the 
Lusztig operators. \end{Rem}
	
	\medskip
	
\begin{Lem}Set \begin{equation} Z_1=Z_{11},
Z_2=Z_{21},Z_3=Z_{12},\textrm{ and }Z_4=Z_{22}.
\end{equation} Then $\forall a\in{\mathbb N}$:
\begin{eqnarray}\label{2cross}Z_4^aZ_1&=&Z_1Z_4^a+(q-q^{-1})
q^{a-1}[a]_qZ_2Z_3Z_4^{a-1}\textrm{ and }\\
Z_4Z_1^a&=&Z_1^aZ_4+(q-q^{-1})
q^{a-1}[a]_qZ_1^{a-1}Z_2Z_3.\end{eqnarray} \end{Lem}

\proof This follows easily by induction from Prop~\ref{qrel}. \qed

	  \bigskip

\section{	${\mathcal A}_q^\pm$ as ${\mathcal 
U}_q({\mathfrak k}^{\mathbb C})$ modules} 
	
\subsection{The quantized Weyl algebra} We wish to describe the natural 
left actions of ${\mathcal 
U}_q({\mathfrak k}^{\mathbb C})$ in these spaces in terms of some 
simple operators given by their matrix representation with respect to a given 
PBW basis. Specifically, introduce, for all 
$i,j\in\{1,2,\dots,n\}$:

\begin{Def} 
\begin{eqnarray}M^o_{ij}(Z_{11}^{a_{11}}\cdot \cdot
Z_{ij}^{a_{ij}}\cdot\cdots
Z_{nn}^{a_{nn}})&=&Z_{11}^{a_{11}}\cdot \cdot
Z_{ij}^{a_{ij}+1}\cdot\cdots Z_{nn}^{a_{nn}}\\
D^o_{ij}(Z_{11}^{a_{11}}\cdot \cdot
Z_{ij}^{a_{ij}}\cdot\cdots
Z_{nn}^{a_{nn}})&=&[a_{ij}]Z_{11}^{a_{11}}\cdot \cdot
Z_{ij}^{a_{ij}-1}\cdot\cdots Z_{nn}^{a_{nn}}\\
H^o_{ij}(Z_{11}^{a_{11}}\cdot \cdot
Z_{ij}^{a_{ij}}\cdot\cdots
Z_{nn}^{a_{nn}})&=&q^{a_{ij}}Z_{11}^{a_{11}}\cdot \cdot
Z_{ij}^{a_{ij}+1}\cdot\cdots Z_{nn}^{a_{nn}}\\
M^o_{ij}(W_{11}^{a_{11}}\cdot \cdot
W_{ij}^{a_{ij}}\cdot\cdots
W_{nn}^{a_{nn}})&=&W_{11}^{a_{11}}\cdot \cdot
W_{ij}^{a_{ij}+1}\cdot\cdots W_{nn}^{a_{nn}}\\
D^o_{ij}(W_{11}^{a_{11}}\cdot \cdot
W_{ij}^{a_{ij}}\cdot\cdots
W_{nn}^{a_{nn}})&=&[a_{ij}]W_{11}^{a_{11}}\cdot \cdot
W_{ij}^{a_{ij}-1}\cdot\cdots W_{nn}^{a_{nn}}\\
H^o_{ij}(W_{11}^{a_{11}}\cdot \cdot
W_{ij}^{a_{ij}}\cdot\cdots
W_{nn}^{a_{nn}})&=&q^{a_{ij}}W_{11}^{a_{11}}\cdot \cdot
W_{ij}^{a_{ij}}\cdot\cdots W_{nn}^{a_{nn}}.\label{32}
\end{eqnarray}\end{Def}
(Notice in particular (\ref{32}).)

\begin{Lem}We have the following formulas for
 all $i,j\in\{1,2,\dots,n\}$:
 \begin{eqnarray}\label{f1}
D^o_{ij}M^o_{ij}-qM^o_{ij}D^o_{ij}&=&(H^o_{ij})^{-1}\\\label{f2}
D^o_{ij}M^o_{ij}-q^{-1}M_{ij}^oD_{ij}^o&
=&H^o_{ij}\\\label{f3} H^o_{ij}D^o_{ij}&=&q^{ -1}
D_{ij}^oH^o_{ij}\\ H^o_{ij}M_{ij}\label{f4}
^o&=&qM_{ij}^oH^o_{ij} .\end{eqnarray}Operators belonging to different nodes 
commute. \end{Lem} \medskip

\begin{Def}We define ${\mathcal W}eyl_q(n,n)$ to be the
algebra generated by the operators $M^o_{ij}, D^o_{ij}$
for all $(i,j)\in \{1,\dots,n\}^2$: \begin{equation}{\mathcal
W}eyl_q(n,n):={\mathbb C}[M^o_{ij}, D^o_{ij};
i,j=1,\dots,n]={\mathcal W}eyl_q(1,1)^{\times n^2}.\end{equation}
Observe that $(H^o_{ij})^{\pm1}\in {\mathcal W}eyl_q(n,n)$
for all $i,j\in\{1,2,\dots,n\}$. \end{Def}

\begin{Rem} This algebra was first studied bu T. Hiashi (\cite{hi}). Recently 
it was studied again in (\cite{kilu}).
\end{Rem}

\subsection{The left actions of ${\mathcal
U}_q({\mathfrak k})$ on ${\mathcal A}_q^{\pm}$}.

We denote the natural left action of a root vector
$X_\mu$  on ${\mathcal
A}_q^{\pm}$, with $X=E,F$, and $\mu\in\Pi$,  by $X^{\pm}_\mu$. We see this 
action as taking place in ${\mathcal U}_q^\pm$ and for this reason introduce 
right multiplication operators $M_{X_\mu}$: ${\mathcal A}_q^\pm\ni u\rightarrow 
uX_{\mu}\in {\mathcal U}_q^\pm$, and analogous operators $M_{K_\mu}$.

\begin{Prop}\label{231}$\textrm{ In }{\mathcal
A}_q^-$:

\begin{eqnarray}\nonumber E^-_{\mu_{k}}&=&\sum_j(-q)
M^o_{k,j}D^o_{k+1,j}H^o_{k,j+1}\dots
H^o_{k,n}(H^o_{k+1,j+1})^{-1}\dots
(H^o_{k+1,n})^{-1}M_{K_{\mu_k}}\ \\&+&
M_{E_{\mu_{k}}}\\\nonumber
F^-_{\mu_{k}}&=&\sum_j
(-q^{-1})M^o_{k+1,j}D^o_{k,j}H^o_{k+1,1}\dots
H^o_{k+1,j-1}(H^o_{k,1})^{-1}\dots (H^o_{k,j-1})^{-1}\\&+&\prod_j
(H^o_{k,j})^{-1}H^o_{k+1,j}M_{F_{\mu_k}}. \end{eqnarray}
Furthermore,\begin{eqnarray}\textrm{ In }{\mathcal
A}_q^-: K_{\mu_k}&=& \prod_j
H^o_{k,j}(H^o_{k+1,j})^{-1}.\label{on-w}\\\textrm{ In
}{\mathcal A}_q^+: K_{\mu_k}&=& \prod_j
H^o_{k+1,j}(H^o_{k,j})^{-1}\label{on-z}.\end{eqnarray}\end{Prop}
	 	  		  
\proof We have the equations \begin{eqnarray}
&E_{\mu_k} T_{\mu_{j}}T_{\mu_{j-1}}\dots
T_{\mu_{1}}(F_\beta)=\\&\left\{\begin{array}{l}
T_{\mu_{j}}T_{\mu_{j-1}}\dots
T_{\mu_{1}}(F_\beta)E_{\mu_k}\textrm{ if }k\neq
j\\(-q)T_{\mu_{j-1}}T_{\mu_{j-2}}\dots
T_{\mu_{1}}(F_\beta)K_{\mu_k}+T_{\mu_{j}}T_{\mu_{j-1}}\dots
T_{\mu_{1}}(F_\beta)E_{\mu_k}\textrm{ if
}k=j.\end{array}\right. \end{eqnarray} and
\begin{eqnarray} &F_{\mu_k}
T_{\mu_{j}}T_{\mu_{j-1}}\dots
T_{\mu_{1}}(F_\beta)=\\&\left\{\begin{array}{l}T_{\mu_{j}}T_{\mu_{j-1}}\dots
T_{\mu_{1}}(F_\beta)F_{\mu_k}\textrm{ if }k\neq j,
j+1\\-q^{-1}T_{\mu_{j+1}}T_{\mu_{j}}\dots
T_{\mu_{1}}(F_\beta)+q^{-1}
T_{\mu_{j}}T_{\mu_{j-1}}\dots
T_{\mu_{1}}(F_\beta)F_{\mu_k}\textrm{ if }k=j+1\\
qT_{\mu_{j}}T_{\mu_{j-1}}\dots
T_{\mu_{1}}(F_\beta)F_{\mu_j}\textrm{ if
}k=j.\end{array}\right. \nonumber\end{eqnarray} The first
equation is clear if $k>j$ and if $k<j$ we can move
$E_{\mu_k}$ past those $T_\ell$ for which $\ell>k$ and
this reduces easily the case to $k=j$. Here it follows
directly from the defining relations of the quantum
group.

The other equation follows analogously. After that we
get

\begin{eqnarray} E_{\mu_k}
W_{i,j}&=&W_{i,j}E_{\mu_k}\textrm{ if }k\neq i-1\\\label{79}
E_{\mu_k}W_{i,j}^a&=&(-q)[a]W_{i-1,j}W_{i,j}^{a-1}K_{\mu_k}+W^a_{i,j}E_{\mu_k}
\textrm{ if }k=i-1\\ F_{\mu_k}
W_{i,j}&=&W_{i,j}F_{\mu_k}\textrm{ if }k\neq i,i-1\\
F_{\mu_k}W_{i,j}^a&=&-q^{-1}[a]W_{i,j}^{a-1}W_{i+1,j}+q^{-a}W^a_{i,j}F_{\mu_k}
\textrm{ if }k=i\\ F_{\mu_k}
W_{i,j}&=&qW_{i,j}F_{\mu_k}\textrm{ if }k= i-1.
\end{eqnarray}

This leads to

\begin{eqnarray}
&E_{\mu_{i-1}}W_{i-1,1}^{a_{i-1,1}}W_{i-1,2}^{a_{i-1,2}}\dots
W_{i-1,j}^{a_{i-1,j}}\dots W_{i-1,n}^{a_{i-1,n}}\cdot
W_{i,1}^{a_{i,1}}W_{i,2}^{a_{i,2}}\dots
W_{i,j}^{a_{i,j}}\dots\\\nonumber&\dots W_{i,n}^{a_{i,n}}\cdot
W_{i+1,1}^{a_{i+1,1}}W_{i+1,2}^{a_{i+1,2}}\dots
W_{i+1,j}^{a_{i+1,j}}\dots
W_{i+1,n}^{a_{i+1,n}}=\\\nonumber&\sum_{j=1}^n
W_{i-1,1}^{a_{i-1,1}}W_{i-1,2}^{a_{i-1,2}}\dots
W_{i-1,j}^{a_{i-1,j}+1}\dots W_{i-1,n}^{a_{i-1,n}}\cdot
W_{i,1}^{a_{i,1}}W_{i,2}^{a_{i,2}}\dots
W_{i,j}^{a_{i,j}-1}\\\nonumber&\dots W_{i,n}^{a_{i,n}}\cdot
W_{i+1,1}^{a_{i+1,1}}W_{i+1,2}^{a_{i+1,2}}\dots
W_{i+1,j}^{a_{i+1,j}}\dots\\\nonumber&\dots
W_{i+1,n}^{a_{i+1,n}}\cdot
(-q)K_{\mu_{i-1}}[a_{i,j}]q^{a_{i-1,j+1}+\dots+a_{i-1,n}+a_{i,j+1}-a_{i,j+1}
-\dots -a_{i,n}}\\\nonumber
&+W_{i-1,1}^{a_{i-1,1}}W_{i-1,2}^{a_{i-1,2}}\dots
W_{i-1,j}^{a_{i-1,j}}\dots W_{i-1,n}^{a_{i-1,n}}\cdot
W_{i,1}^{a_{i,1}}W_{i,2}^{a_{i,2}}\dots
W_{i,j}^{a_{i,j}}\\\nonumber&\dots W_{i,n}^{a_{i,n}}\cdot
W_{i+1,1}^{a_{i+1,1}}W_{i+1,2}^{a_{i+1,2}}\dots
W_{i+1,j}^{a_{i+1,j}}\dots
W_{i+1,n}^{a_{i+1,n}}E_{\mu_{i-1}}. \end{eqnarray}
Notice that the $q$ exponents arise from the rearranging
of terms (e.g. $W_{i-1,j}$) and $K_{\mu_k}$.

There are analogous considerations for $F_{\mu_k}$. \qed

\bigskip

The case of the $Z_{ij}$s is similar: \smallskip First
observe the very useful formula
\begin{equation}ad(E_{\mu_{1}})ad(E_{\mu_{2}})ad(E_{\mu_{1}})(E_\beta)=0.
\end{equation}

This follows from a lengthy computation based on the
Serre relations. Then observe that

\begin{eqnarray} \label{80}&E_{\mu_k}
ad(E_{\mu_{j}})ad(E_{\mu_{j-1}})\dots
ad(E_{\mu_{1}})(E_\beta)\\&=\left\{\begin{array}{l}ad(E_{\mu_{j}})ad(E_{\mu_{j-1
}})\dots
ad(E_{\mu_{1}})(E_\beta)E_{\mu_k} \textrm{ if }k\neq j,
j+1\\ \ad(E_{\mu_{j+1}})ad(E_{\mu_{j}})\dots
ad(E_{\mu_{1}}) (E_\beta)\\+q^{-1} ad(E_{\mu_{j}})\dots
ad(E_{\mu_{1}})(E_\beta)E_{\mu_k} \textrm{ if }k=j+1\\\nonumber
q\cdot ad(E_{\mu_{j}})ad(E_{\mu_{j-1}})\dots
ad(E_{\mu_{1}})(E_\beta)E_{\mu_j}\textrm{ if
}k=j.\end{array}\right. \end{eqnarray}

\begin{eqnarray} &F_{\mu_k}
ad(E_{\mu_{j}})ad(E_{\mu_{j-1}})\dots
ad(E_{\mu_{1}})(E_\beta)\\&=\left\{\begin{array}{l}ad(E_{\mu_{j}})ad(E_{\mu_{j-1
}})\dots
ad(E_{\mu_{1}})(E_\beta)F_{\mu_k} \textrm{ if }k\neq j\\\nonumber
ad(E_{\mu_{j-1}})\dots ad(E_{\mu_{1}})
(E_\beta)K^{-1}_{\mu_k}+ad(E_{\mu_{j}}) \dots
ad(E_{\mu_{1}})(E_\beta)F_{\mu_k} \textrm{ if
}k=j.\end{array}\right. \end{eqnarray}

The form of the operators follow easily from this, in a
way similar to the case of ${\mathcal A}_q^-$.

\smallskip

We have then proved

\begin{Prop} \label{k-prod}\begin{eqnarray}\label{232} &E^+_{\mu_k}=\\&\sum_j
D^o_{k,j}M^o_{k+1,j}(H^o_{k,1})^{-1}\dots\nonumber
(H^o_{k,j-1})^{-1}H^o_{k+1,1}\dots H^o_{k+1,j-1}\\\nonumber&{ +
(H^o_{k,1})^{-1}\dots (H^o_{k,n})^{-1}H^o_{k+1,1}\dots H^o_{k+1,n}
M_{E_{\mu_k}}},\\
&F^+_{\mu_k}=\\&\sum_j
D^o_{k+1,j}M^o_{k,j}(H^o_{k+1,j+1})^{-1}\dots
(H^o_{k+1,n})^{-1}H^o_{k,j+1}\dots H^o_{k,n}
{\cdot
M_{K^{-1}_{\mu_k}}}.\end{eqnarray} \end{Prop}

\medskip

We use the antipode $S=S_J$ from (\cite{jan}, p.34). Specifically
\begin{equation}
\forall \alpha\in\Pi: S_J(E_\alpha)=-K_\alpha^{-1}E_\alpha,\  
S_J(F_\alpha)=-F_\alpha K_\alpha, S_J(K_\alpha)=K_\alpha^{-1}.
\end{equation}

\smallskip

We shall later study dual modules. Here we recall the 
definition\begin{equation} (u^T
v',v):=(v',S(u)v).\end{equation}

\smallskip

With this in mind, we observe that we have, modulo right actions by $E_\mu$ and 
$F_\mu$, \begin{Cor}
\begin{eqnarray}&(-K_{\mu_k}^{-1})E^-_{\mu_k}=\\\nonumber&\sum_j
q
(H^o_{k,j})^{-1}M^o_{k,j}H^o_{k+1,j}
D^o_{k+1,j}H^o_{k+1,1}\dots
H^o_{k+1,j-1}(H^o_{k,1})^{-1}\dots (H^o_{k,j-1})^{-1}\label{233}\\
&q\cdot F^-_{\mu_{k}}(-K_{\mu_k})=\\\nonumber&\sum_j
M^o_{k+1,j}(H^o_{k+1,j})^{-1}D^o_{k,j}H^o_{k,j}\cdot
(H^{o}_{k+1,j+1})^{-1}\dots (H^{o}_{k+1,n})^{-1}H^o_{k,j+1}\dots
H_{k,n}.\end{eqnarray} \end{Cor}

Of course, there are similar formulas for
the actions of $E^\pm_\nu$ and $F^\pm_\nu$. We omit those as they are entirely 
similar, easily deducible, and since, for our purposes, they do not add 
anything new.

	  \bigskip

	  \section{Duality}

\subsection{The $q$ Killing form d'apr\`es Jantzen}
The $q$ version of the Killing form was introduced by M. Rosso (\cite{rosso}). 
Here we follow the comprehensive study offered in (\cite{jan}); a study that 
relies on the approaches offered in (\cite{tani}) and (\cite{jo-le}).
	  
We cite:

\begin{Prop}[\cite{jan} Proposition~ 6.12] \label{p612}There exists a unique
bilinear pairing $(\cdot,\cdot)_{J}$: ${\mathcal
U}^{\geq0}\times {\mathcal U}^{\leq0}\rightarrow k$ such
that for all $x,x'\in {\mathcal U}^{\geq0}$, all
$y,y'\in {\mathcal U}^{\leq0}$, all $\mu,\nu\in{\mathbf
Z}\Psi$, and all $\alpha,\beta\in\Pi$, \begin{eqnarray}
(y,xx')_{J}=(\triangle(y),x'\otimes
x)_{J}&,&(yy',x)_{{J}}=(y\otimes y',\triangle(x))_{{J}},\\
(K_\mu,K_\nu)_{{J}}=q^{-(\mu,\nu)}&,&(F_\alpha,E_\beta)_{{J}}=-\delta_{\alpha,
\beta}
(q_\alpha-q_\alpha^{-1})^{-1}\\
(K_\mu,E_\alpha)_{{J}}=0&,&
(F_\alpha,K_\mu)_{{J}}=0.
\end{eqnarray} (The form is extended to tensor
products in the natural way). \end{Prop}

\begin{Rem} The form is unique {\em for the given}
$\triangle$. \end{Rem}

\medskip

The following follows easily:
\begin{Cor}\label{K}For all $x\in {\mathcal U}^{>0}$, all
$y\in {\mathcal U}^{<0}$, all $\mu,\nu\in{\mathbf
Z}\Psi$,
\begin{equation}(yK_\mu,x K_\nu)_{J}=(y,x)_Jq^{-(\mu,\nu)}.
\end{equation}
\end{Cor}

\medskip

\begin{Prop}[\cite{jan} Corollary in \S8.30] If $q$ is not a root of unity, the 
form restricted to ${\mathcal U}_q^{-\mu}\times {\mathcal U}_q^{\mu}$ is 
non-degenerate for any $\mu\in{\mathcal L}^+$.
\end{Prop}

The following also holds. It can be proved using the same argument as in 
(\cite[Proposition~6.21]{jan}).

\begin{Prop}If $q$ is not a root of unity,
\begin{equation}(K_\alpha , 
K_\gamma)_J=q^{-(\alpha,\gamma)}\end{equation}extends to a non-degenerate 
bilinear form on ${\mathcal U}_q^0\times {\mathcal U}_q^0$.\end{Prop}

\begin{Prop}[\cite{jan} 6.13 (3)]
\begin{equation}(Ad(K_\alpha)x,y)_{J}=(x, Ad(K_\alpha^{-1}y)_{J}.
\end{equation}
\end{Prop}

\medskip

\begin{Prop}[\cite{jan} Proposition~8.29]Let
$\omega=s_{\alpha_1}s_{\alpha_2}\dots s_{\alpha_t}$ be a
reduced representation. Then
$T_{\alpha_1}(F_{\alpha_2}^b)=(T_{\alpha_1}(F_{\alpha_2}))^b$.
Furthermore, \begin{equation}\nonumber\left(T_{\alpha_1}T_{\alpha_2}\dots
T_{\alpha_t}(F_{\alpha_t}^{b_t})\dots
T_{\alpha_1}(F_{\alpha_2}^{b_2})F_{\alpha_1}^{b_1},T_{\alpha_1}T_{\alpha_2}\dots
T_{\alpha_t}(E_{\alpha_t}^{a_t})\dots
T_{\alpha_1}(E_{\alpha_2}^{a_2})E_{\alpha_1}^{a_1}\right)_{J}\end{equation}
\begin{equation}=\prod_i\delta_{a_i,b_i
}\prod_{j=1}^t
(F_{\alpha_j}^{a_j},E_{\alpha_j}^{a_j})_{J}.\end{equation}\end{Prop}

\begin{Lem}[\cite{jan} (4) p. 114]
\begin{equation}\left(F_{\alpha_i}^{a_i},E_{\alpha_i}^{a_i}\right)_{J}
=(-1)^{a_i}q_{\alpha_i}^{a_i(a_i-1)/2}[a_i]_{\alpha_i}!/(q_{\alpha_i}-q_{
\alpha_i}^{-1})^{\alpha_i}\label{93}.\end{equation} \end{Lem} \medskip

According to our definitions, the following is immediate:

\begin{Cor}
\begin{equation}\left(W_{11}^{b_{11}}W_{12}^{b_{12}}\cdots
W_{1n}^{b_{1n}}W_{21}^{b_{21}}\cdots
W_{nn}^{b_{nn}},Z_{11}^{a_{11}}Z_{12}^{a_{12}}\cdots
Z_{1n}^{a_{1n}}Z_{21}^{a_{21}}\cdots
Z_{nn}^{a_{nn}}\right)_{J}\nonumber\end{equation}\begin{equation}=\prod_{(ij),
(rs)
}\delta_{(ij),(rs)}\prod_{(ij)}
\left(\frac{-1}{q-q^{-1}}\right)^{a_{ij}}q^{a_i(a_i-1)/2}[a_i]!\end{equation}
\label{jan-form}
\end{Cor}

\bigskip

 \subsection{More formulas from Jantzen}

\medskip

\begin{Def}[\cite{jan} \S6.14] We here define the important operators
$r_\alpha$ and $r'_\alpha$ in ${\mathcal U}_q^+$: \begin{eqnarray}x\in
U_\mu^+:\Delta(x)&=&x\otimes 1+
\sum_{\alpha\in\Pi}r_\alpha(x)K_\alpha\otimes
E_\alpha+(rest).\\ x\in U_\mu^+:\Delta(x)&=&K_\mu\otimes
x+ \sum_{\alpha\in\Pi}E_\alpha K_{\mu-\alpha}\otimes
r'_\alpha(x)+(rest).\end{eqnarray} \end{Def}

\begin{Prop}[\cite{jan} \S6.14] If $x\in{\mathcal U}_{\mu}$ then
$r_\alpha(x),r'_\alpha(x)\in{\mathcal U}_{-\alpha+\mu}$.
Let furthermore $x'\in{\mathcal U}_{\mu'}$. Then
\begin{eqnarray}r_\alpha(1)=r'_\alpha(1)=0&;&
r_\alpha(E_\beta)=r'_\alpha(E_\beta)=\delta_{\alpha,\beta},
\\r'_\alpha(xx')&=&r'_\alpha(x)x'+q^{ (\alpha,\mu)}
xr'_\alpha(x'),\\
r_\alpha(xx')&=&q^{(\alpha,\mu')}r_\alpha(x)x'+xr_\alpha(x'),\\\label{65}
(u^-F_\alpha,u^+)_{J}&=&(F_\alpha,E_\alpha)_{J}(u^-,
r_\alpha(u^+))_{J},\textrm{ and}\\\label{66} (F_\alpha
u^-,u^+)_{J}&=&(F_\alpha,E_\alpha)_{J}(u^-,r'_\alpha(u^+))_{J}.
\end{eqnarray}
\end{Prop}

There are analogous operators $r_\alpha$ and $r'_\alpha$
in ${\mathcal U}_q^-$:

\begin{Def}[\cite{jan} \S6.15] 
\begin{eqnarray}y\in U_{-\mu}^-:\Delta(y)&=&y\otimes
K_{\mu}^{-1}+ \sum_{\alpha\in\Pi}r_\alpha(y)\otimes
F_\alpha K_{\mu-\alpha}^{-1}+(rest).\\ y\in
U_{-\mu}^-:\Delta(y)&=&1\otimes y+
\sum_{\alpha\in\Pi}F_\alpha \otimes
r'_\alpha(y)K_{\alpha}^{-1}+(rest).\end{eqnarray}
\end{Def}

\medskip

\begin{Prop}[\cite{jan} \S6.15] If $y\in{\mathcal U}^-_{-\mu}$ then
$r_\alpha(y),r'_\alpha(y)\in{\mathcal
U}^-_{\alpha-\mu}$. Let furthermore $y'\in{\mathcal
U}^-_{-\mu'}$. Then
\begin{eqnarray}r_\alpha(1)=r'_\alpha(1)=0&;&
r_\alpha(F_\beta)=r'_\alpha(F_\beta)=\delta_{\alpha,\beta},
\\r_\alpha(yy')&=&r_\alpha(y)y'+ q^{ (\alpha , \mu) }
yr_\alpha(y'),\\
r'_\alpha(yy')&=&q^{(\alpha,\mu')}r'_\alpha(y)y'+yr'_\alpha(y'),\\
\forall \alpha,\beta\in\Pi:r_\alpha\circ
r'_\beta&=&r'_\beta\circ r_\alpha \ (\textrm{Jantzen p. 218 l. 1}),\\\label{58}
(y^-,E_\alpha
y^+)_{J}&=&(F_\alpha,E_\alpha)_{J}(r_\alpha(y^-), y^+)_{J}\textrm{ 
and}\\\label{59}
(y^-,y^+E_\alpha
)_{J}&=&(F_\alpha,E_\alpha)_{J}(r'_\alpha(y^-), y^+)_{J}.\label{60}
\end{eqnarray}
\end{Prop}

\smallskip

Observe that
\begin{equation}r_\nu(T_\nu(E_\beta))=0=r'_\nu(T_\nu(F_\beta)).\end{equation}
\bigskip

From this we easily get the following special cases:

\begin{Lem}
\begin{equation}\begin{array}{lclclclclcl} r_\beta
W_{ij}&=&\delta_{1,i}\delta_{1,j},&&
r_\beta^nW_{11}^n&=&[[n]]_q!,\\
r_{\mu_i}W_{k,j}&=&\delta_{k,i+1}(-\gamma)W_{i,j},&&
r_{\mu_i}&=&0 \textrm{ on }{\mathcal A}_q^+,\\ r'_{\mu_i}&=&0 \textrm{
on }{\mathcal A}_q^-,&& r'_\beta Z_{ij}&=&\delta_{1,i}\delta_{1,j},\\
{r'_\beta}^nZ_{11}^n&=&[[n]]_q!,\qquad\textrm{ and}&&
r'_{\mu_i}Z_{k,j}&=&\delta_{k,i+1}(q^{-1}\gamma)Z_{i,j}.
\end{array}\end{equation}
\end{Lem}

Jantzen's formulas in Corollary~\ref{jan-form} follows readily from these.

\bigskip

\begin{Rem}The bilinear form $(\cdot,\cdot)_{J}$ has a singularity at $q=1$. We 
will renormalize it later in a fixed PBW basis, but for now we keep it because 
of (\ref{58}) and (\ref{59}). Towards the end of his book, in \S10.16, Jantzen 
introduces a renormalized bilinear form. This we will not use, since it, for 
our purposes, is more difficult to use.
\end{Rem}

\bigskip

\subsection{The left action of $E_\alpha$ in ${\mathcal
U}_q^-$.}

\begin{Prop}[\cite{jan} Lemma~6.17] With $\alpha\in\Pi$,
$\mu\in{\mathbb Z}\Phi, y\in{\mathcal U}^-_{-\mu}$, and
$x\in{\mathcal U}^+_{\mu}$: \begin{eqnarray}E_\alpha
y-y
E_\alpha&=&\frac1{q_\alpha-q_\alpha^{-1}}\left(K_\alpha
r_\alpha(y)-r'_\alpha(y)K_\alpha^{-1}\right)\\ x
F_\alpha-F_\alpha
x&=&\frac1{q_\alpha-q_\alpha^{-1}}\left(r_\alpha(x)K_\alpha
-K_\alpha^{-1}r'_\alpha(x)\right). \end{eqnarray}
\end{Prop}(Jantzen proves this by a simple induction
argument.) 

\smallskip

From (\ref{93}), (\ref{59}), and (\ref{60}) one easily deduces:
\begin{Cor}for elements $u^-_{-\mu}\in{\mathcal U}^-_{-\mu}$ and 
$u^+_{\mu-\alpha}\in{\mathcal U}^+_{\mu-\alpha}$ we have
\begin{equation}
\left(E_\alpha u^-_{-\mu}-u^-_{-\mu}E_\alpha, 
u^+_{\mu-\alpha}\right)_{J}=-\left(u^-_{-\mu}, (E_\alpha 
K_\alpha^{-1})u^+_{\mu-\alpha}-u^+_{\mu-\alpha}(K_\alpha 
E_\alpha)\right)_{J}.\label{116}\end{equation}
\end{Cor}

	\bigskip

\section{Duality reconsidered,   especially the  ${\mathcal U_q}({\mathfrak 
k})$ modules ${\mathcal A}_q^\pm$} 
We consider here
non-singular pairings between a (highest weight) module
and some other module. By this we mean in general a complex valued 
non-degenerate bilinear form taking inputs from two modules such that the 
second module is the dual module of the first according to the
pairing. On some level, there is of course only one dual module, but they may 
be given in different realizations.

\subsection{More on the Jantzen pairing}

It is easy to see, since $r_\alpha=0$ in ${\mathcal
A}_q^+$ for any simple compact root $\alpha$, that

\begin{Lem}$\forall u^+_z\in {\mathcal A}_q^+,u^-_w\in
{\mathcal A}_q^-, u^+_{\mathfrak k}\in {\mathcal
U}_q^+$, and $u^-_{\mathfrak k}\in {\mathcal U}_q^-$:
\begin{equation}(u^-_{w}u^-_{\mathfrak k}, u^+_{z}u^+_{\mathfrak
	k})_{J}=(u^-_{w}, u^+_{z})_{J}(u^-_{\mathfrak
k},u^+_{\mathfrak k})_{J}.\end{equation} \end{Lem}

\medskip

We will use this version of the pairing because there
are some very simple formulas for the duals of the
operators $r_\alpha,r'_\alpha$ to be studied later. However, there will also be 
modified pairings:

\subsection{Pairings between ${\mathcal A}_q^-$ and
${\mathcal A}_q^+$.}

\begin{Def}\begin{eqnarray}[[a]]_q&=&1+q^2+\cdots + q^{2a-2},\\
 \left[a\right]_q&=&q^{-a+1}+\cdots+q^{a-1},\textrm{ and
}\\\{\{a\}\}_q&=&1+q^{-2}+\cdots + q^{-(2a-2)}.\end{eqnarray} \end{Def}

If
${\mathbf a}=(a_{11},a_{12},\dots,a_{nn})$, ${\mathbf
b}=(b_{11},b_{12},\dots,b_{nn})$, and ${\mathbf
c}=(c_{11},c_{12},\dots,c_{nn})$ we let $[[{\mathbf
c}]]!=\prod_{ij}([[c_{ij}]]!), \delta_{{\mathbf
a},{\mathbf
b}}=\prod_{ij}\delta_{a_{ij},b_{ij}}, \vert{\mathbf
a}\vert=\sum_{ij}a_{ij}$, { and }$
Z^{\mathbf c}=Z_{11}^{c_{11}}Z_{12}^{c_{12}}\cdots
Z_{nn}^{c_{nn}}$.

\medskip

\begin{Def}The following bilinear forms, indexed by $J,K$, and $L$, will be 
considered:\begin{eqnarray}\label{41}
\left(Z^{\mathbf a},W^{\mathbf 
b}\right)_{J}&=&\left(\frac{-1}{q-q^{-1}}\right)^{\vert{\mathbf
a}\vert}
\delta_{{\mathbf a},{\mathbf b}} [[{\mathbf a}]]!,\\\left(Z^{\mathbf
a},W^{\mathbf b}\right)_L&=&\left(\frac{-1}{q-q^{-1}}\right)^{\vert{\mathbf
a}\vert}\delta_{{\mathbf a},
{\mathbf b}}[{\mathbf a}]!, \textrm{ and
}\\
\left(Z^{\mathbf a},W^{\mathbf
b}\right)_K&=&\left(\frac{-1}{q-q^{-1}}\right)^{\vert{\mathbf
a}\vert}\delta_{{\mathbf a},{\mathbf
b}}\{\{{\mathbf a}\}\}.
\label{43}\end{eqnarray} \end{Def}

\medskip

Corresponding to these forms we introduce two more
families of differential operators (we include the old
one for convenience): \begin{Def}$\forall i,j=1,\dots,n$
\begin{eqnarray} {\mathbb
D}^o_{ij}Z_{ij}^{a_{ij}}&=&[[a_{ij}]]Z_{ij}^{a_{ij}-1},\\
D^o_{ij}Z_{ij}^{a_{ij}}&=&[a_{ij}]Z_{ij}^{a_{ij}-1}, \textrm{ and
}\\
{\mathcal
D}^o_{ij}Z_{ij}^{a_{ij}}&=&\{\{a_{ij}\}\}Z_{ij}^{a_{ij}-1}.
\end{eqnarray} Clearly, ${\mathbb
D}^o_{ij}=H^o_{ij}D_{ij}^o$ and ${\mathcal
D}^o_{ij}=(H^{o}_{ij})^{-1}D_{ij}^o$. \end{Def}

\medskip

Of course, we have formulas analogous to
(\ref{f1}-\ref{f4}):

\begin{Lem}
\begin{eqnarray}\label{f5} {\mathbb
D}^o_{ij}M^o_{ij}-q^2M^o_{ij}{\mathbb
D}^o_{ij}&=&I,\\\label{f6} {\mathbb
D}^o_{ij}M^o_{ij}-M_{ij}^o{\mathbb D}_{ij}^o&
=&(H^o_{ij})^2,\\\label{f7} H^o_{ij}{\mathbb D}^o_{ij}&=&q^{
-1} {\mathbb D}_{ij}^oH^o_{ij}, \textrm{ and
}\\ H^o_{ij}M_{ij}\label{f8}
^o&=&qM_{ij}^oH^o_{ij}.\end{eqnarray}
and \begin{eqnarray}\label{f5} {\mathcal
D}^o_{ij}M^o_{ij}-q^{-2}M^o_{ij}{\mathcal
D}^o_{ij}&=&I,\\\label{f6} {\mathcal
D}^o_{ij}M^o_{ij}-M_{ij}^o{\mathcal D}_{ij}^o&
=&(H^{o}_{ij})^{-2},\\\label{f7} H^o_{ij}{\mathcal
D}^o_{ij}&=&q^{ -1} {\mathcal D}_{ij}^oH^o_{ij}, \textrm{ and
}\\
H^o_{ij}M_{ij}\label{f8} ^o&=&qM_{ij}^oH^o_{ij}.\end{eqnarray}
\end{Lem}

\bigskip

\bigskip

\subsection{Returning to the  ${\mathcal U}_q({\mathfrak k}^{\mathbb C})$ 
modules ${\mathcal A}_q^\pm$.}

\label{5.3}

We denote by $\widehat{T}^X$ the dual operator of $T$ according
to some pairing $(\cdot,\cdot)_X$. Specifically,
\begin{eqnarray}\forall u,v\in{\mathcal A}^+_q\times {\mathcal A}^-_q: 
(\widehat{T}^Xv,u)_X&=&(v,Tu)_X, \textrm{ or}\\
(Tv,u)_X&=&(v,\widehat{T}^Xu)_X.
\end{eqnarray}
 
   We assume $\forall i,j\in\{1,2,\dots,n\}: 
\widehat{H^o_{ij}}^X=H^o_{ij}$ and seek ``self duality'' in the sense
of Jantzen p. 120 for the operators from Propositions \ref{231} and \ref{k-prod}.

\medskip

\begin{Prop}\label{5.4} There is the following equivalence:
\begin{equation}\forall \xi\in \Pi_c:\ \left.\begin{array}{ccc}
\widehat{S_J(E_\xi^-)}^X&=&E_\xi^+
\\\widehat{S_J(F_\xi^-)}^X&
=&F_\xi^+\\
\widehat{S_J^{-1}(E_\xi^+)}^X&=&E_\xi^-\\\widehat{S_J^{-1}(F_\xi^+)}
^X&=&F_\xi^- \end{array}\right\}\Leftrightarrow
\forall 
i,j\in\{1,2,\dots,n\}:\widehat{M_ij}^X=\kappa 
H^o_{ij}D^o_{ij}\end{equation}holds,
where $\kappa$ is a complex constant. 
\end{Prop}

\proof We need to compare terms in (\ref{232}) and  (\ref{233}). It follows 
easily, by looking at the operators involving $E^\pm_{\mu_k}$ that, for $y$ 
representing any double index,  $\widehat{ D^o_y}^X=(H^o_y)^{-1}M^o_yq^\alpha$ 
and $\widehat{M^o_y}^X= H^o_yD^o_yq^\beta$ with $\alpha+\beta=1$. It follows 
from these that $T\rightarrow \widehat{T}^X$ is an idempotent. Moreover, the 
equations are equivalent and are also equivalent to those arising from  
$F^\pm_{\mu_k}$.

\medskip

\begin{Rem}J.C. Jantzen introduces in \S6.20 an altered bilinear form 
$<\cdot,\cdot>$ on ${\mathcal U}_q$ in which $ad$ is self dual (his Proposition 
6.20). This form does not restrict to $(\cdot,\cdot)_J$. It is not clear that 
self-duality of restrictions of $\ad$  to be desired. There is also a 
non-uniqueness in the sense that the co-product $\triangle$ may be altered. 
Furthermore, notice that \begin{equation}E\rightarrow E_x=EK^x;\ F\rightarrow
F_x=K^{-x}F\end{equation} defines an automorphism of the quantized enveloping 
algebra (preserves the $q$-Serre relations). \end{Rem}

The following is straightforward:
\begin{Lem}If for some fixed $i,j$ $\widehat{W_{ij}}^X=\kappa H^o_{ij}D^o_{ij}$ 
for some $\kappa\neq0$ and some pairing $X$ then
\begin{equation}\forall a,b\in{\mathbb N}_0:(W_{ij}^a,Z_{ij}^b)_X=
\delta_{a,b}\kappa^a[[a]]!\end{equation}
In other words, the form above is essentially the form $(\cdot,\cdot)_{J}$.
\end{Lem}

	  \medskip
	  
	  Notice that $\widehat{M}^L=D$. We shall see in later sections that 
there are difficulties with the form $(\cdot,\cdot)_{J}$, whereas the two other 
forms indexed by $L$ and $K$, respectively, behave very nicely.

\bigskip

\begin{Rem}
We will also later define differential operators as duals of either left or 
right multiplication operators. Here we take the stance of using ``flat'' 
dualities without using antipodes. This is so because there is no natural 
antipode for, especially, right multiplication operators.
\end{Rem}

\subsection{Change of basis}

Let ${\mathbb
A}W^\alpha=q^{-\frac{\alpha(\alpha-1)}2}W^\alpha$.

$\left(Z^\alpha,W^\beta\right)_L=\left(Z^\alpha,{\mathbb
A}W^\beta\right)_{J}$. Clearly, $\widehat{{\mathbb
A}}^{J}={\mathbb A}$.

 Then
\begin{equation}\widehat{T}^L=\widehat{\left({\mathbb A}T{\mathbb
A}^{-1}\right)}^{J} = {\mathbb A}^{-1}\widehat{T}^{J}{\mathbb
A} .\end{equation}

 Suppose we are given  $\widehat{T}^J$ in terms of sums of monomials in the 
operators $L,M,H$, then, 
to get $\widehat{T}^L$, we just need to make the replacements, in the given 
expression,

\begin{eqnarray}\label{g1} M\textrm{ by }\left({\mathbb
A}^{-1}M{\mathbb A}\right)&=&MH^o,\\ \textrm{ and }D \textrm{ by } 
\left({\mathbb
A}^{-1}D{\mathbb A}\right)&=&(H^o)^{-1}D={\mathcal
D}.\label{g3}
 \end{eqnarray}

Of course, \begin{eqnarray}\widehat{T}^{K}=\widehat{\left({\mathbb A}^2T{\mathbb
A}^{-2}\right)}^{J}&=& {\mathbb A}^{-2}\widehat{T}^{J}{\mathbb
A}^2 ,\textrm{ and}\\ \label{h1} \left({\mathbb
A}^{-2}M{\mathbb A}^2\right)&=&M(H^o)^2,\\ \left({\mathbb
A}^{-2}D{\mathbb A}^2\right)&=&(H^o)^{-2}D.\label{h3}
\end{eqnarray}

\medskip

\begin{Lem} The prescriptions (\ref{g1} - \ref{g3})
above extend to an automorphism of ${\mathcal W}eyl_q(n,n)$.
This is the change-of-basis automorphism. \end{Lem}

\pof This follows immediately from (\ref{f1} -
\ref{f4}). \qed

\medskip

\begin{Rem}
The change-of-basis map ${\mathbb A}$ on ${\mathcal A}^+_q$ can also be viewed 
as a change-of-variable transformation, albeit in an enlarged algebra: 
Let\begin{equation}\widetilde{\mathcal A}^+_q ={\mathcal A}\times_s{\mathbb 
C}[H_{ij}^{\pm1}\mid 
ij=1,\dots,n].\end{equation}Then\begin{equation}(q^{\frac{1}{2n}}Z_{ij}
(H^o_{ij})^{-1})^{a_{ij}}
=q^{-\frac{a_{ij}(a_{ij}-1)}2}(Z_{ij})^{a_{ij}}q^{\frac{a_{ij}}{2n}}(H^o_{ij})^
{-a_{ij}}.\end{equation}In a scalar module (where $\Lambda(\xi)=0$ for all 
$\xi\in \Pi_c$ and $\Lambda(\beta)=\lambda\in{\mathbb R}$) one may set 
$H_{ij}v_\Lambda=q^{\frac{1}{2n}}v_\Lambda$ to enlarge the representation to 
these elements. Of course, further modifications will have to be added to our 
algebras for this 
to make sense - all operators of the form 
$D_{i,j}M^o_{i+1,j}$ should also be included. 

This idea will be persued elsewhere.
\end{Rem}

\bigskip

\section{The full pairing}

We are interested in studying duals of generalized Verma
modules. 

\subsection{Highest Weights}
A highest weight vector of a module of ${\mathcal
U}_q({\mathfrak k}^{\mathbb C})$ has a highest weight
$\Lambda$ and a highest weight vector $v_\Lambda\neq{\bf 0}$ for
which \begin{equation}\begin{array}{cccccccc} \forall
i=1,\dots,n-1:
K_{\mu_i}^{\pm1}=q^{\pm\lambda^\mu_i}v_\Lambda,&
 K_{\nu_i}^{\pm1}=
q^{\pm \lambda^\nu_i}v_\Lambda,&\textrm{ and } 
K_{\beta}^{\pm1}=q^{\pm\lambda}v_\Lambda.
\end{array}\end{equation}
We set 
$\Lambda=((\lambda^\mu_1,\dots,\lambda^\mu_{n-1}),(\lambda^\nu_1,\dots,
\lambda^\nu_{n-1});\lambda)=(\Lambda_L,
\Lambda_R,\lambda)$. We assume:
\begin{equation}\label{inte}\forall i=1,\dots,n-1:\ 
\lambda^\nu_i,\lambda^\mu_j\in {\mathbb N}_0.
\end{equation}

As a vector space, $V_\Lambda=V_{\Lambda_L}\otimes V_{\Lambda_R}$ where 
$V_{\Lambda_L}$ and $V_{\Lambda_R}$ are highest weight representations of
${\mathcal
U}_q({\mathfrak k}_L^{\mathbb C})$ and ${\mathcal
U}_q({\mathfrak k}_R^{\mathbb C})$, respectively, of highest weights 
$\Lambda_L=(\lambda^\mu_1,\dots,\lambda^\mu_{n-1})$ and 
$\Lambda_R=(\lambda^\nu_1,\dots,\lambda^\nu_{n-1})$, 
respectively. The highest weight vector can then be written as 
$v_\Lambda=v_{\Lambda_L}\otimes v_{\Lambda_R}$ with the stipulation that 
$K_\beta^{\pm1}v_{\Lambda_L}\otimes v_{\Lambda_R}=q^{\pm\lambda}v_{\Lambda_L}
\otimes v_{\Lambda_R}$.

\smallskip

The condition (\ref{inte}) is an integrality condition and gives rise to a 
finite-dimensional ${\mathcal U}_q({\mathfrak k})$ module $V_\Lambda$ in the 
following usual way:

${\mathcal I}_L({\mathfrak k})$ denote the left  
${\mathcal U}_q({\mathfrak k})$ invariant subspace in ${\mathcal
U}_q({\mathfrak k}^-)v_\Lambda$ generated by
$\{(F_\gamma)^{\Lambda(\gamma)+1}v_\Lambda\}$ and 
let\begin{equation}V_\Lambda={\mathcal U}_q({\mathfrak k})v_\Lambda/{\mathcal
I}_L({\mathfrak k}).\end{equation}

The dual module $V'_\Lambda$ is the highest weight module 
$V_{-\omega_0(\Lambda)}$. However, it is also convenient to view $V'_\Lambda$ 
as 
a lowest weight module $V^o_{\Lambda'}$ characterized by a non-zero lowest 
weight vector $v^o_{\Lambda'}$ for which
\begin{eqnarray}{\mathcal U }({\mathfrak
k})^-_qv^o_{\Lambda'}&=&0\\
\forall\alpha\in\Pi_c:K_\alpha^{\pm1}
v^o_{\Lambda'}&=&q^{\pm\Lambda'(\alpha)}
v^o_{\Lambda'}\\
\forall\alpha\in\Pi_c:\Lambda'(\alpha)&\in&-{\mathbb
N}_0.\label{44-o} 
\end{eqnarray}

\medskip

The following is elementary

\begin{Lem}\begin{equation}\Lambda'=-\Lambda.\end{equation}
\end{Lem}

\begin{Rem}
One may consider a modified version $(\cdot,\cdot)_{mod}$ of 
$(\cdot,\cdot)_{J}$ 
on ${\mathcal U}_q^-\times {\mathcal U}_q^+$, where 
$\forall u^+_z\in{\mathcal A}_q^+,\forall u^-_w\in{\mathcal A}_q^-,\forall 
u^{\pm}_{\mathfrak k}\in{\mathcal U}_q^{\pm},\forall \lambda,\mu\in 
{\mathfrak L}:$ \begin{equation}
(u^-_wu^-_{\mathfrak k}K_\lambda, 
u^+_zu^+_{\mathfrak 
k}K_\mu)_{mod}=q^{(\Lambda(\lambda)-\Lambda(\mu))}(u^-_wu^-_{\mathfrak k}, 
u^+_zu^+_{\mathfrak k})_{mod}\end{equation} such that the finite-dimensional 
modules 
$V_\Lambda$ and $V^o_{\Lambda'}$ occur naturally as duals in this setting.  
Furthermore, it holds in the same generality as above that  
\begin{equation}(K_\lambda u^-_wu^-_{\mathfrak k}, 
u^+_zu^+_{\mathfrak k})_{mod}=( u^-_wu^-_{\mathfrak k}, K_\lambda^{-1}
u^+_zu^+_{\mathfrak k})_{mod}.\end{equation}
\end{Rem}

\medskip

\subsection{Generalized Verma modules}

Consider a finite dimensional module
$V_\Lambda=V_{\Lambda_L,\Lambda_R,\lambda}$ over
${\mathcal U}_q({\mathfrak k}^{\mathbb C})$ with highest weight is defined by
$\Lambda=(\Lambda_L,\Lambda_R,\lambda)$.

 We extend
such a module to a ${\mathcal U}_q({\mathfrak k}^{\mathbb
C})\cdot {\mathcal A}_q^+$ module, by the same name, by letting
${\mathcal A}_q^+$ act trivially in $V_\Lambda$.

\begin{Def} The quantized generalized Verma module
$M(V_\Lambda)$ is given by
\begin{equation}{\mathcal M}(V_\Lambda)={\mathcal U}_q({\mathfrak
g}^{\mathbb C})\bigotimes_{{\mathcal U}_q({\mathfrak
k}^{\mathbb C}){\mathcal A}_q^+}V_\Lambda \end{equation}
with the natural action from the left. We denote the corresponding 
representation by $L_\Lambda(u)$ for 
$u\in{\mathcal U}_q$.\end{Def} \bigskip

As a vector space, even as a ${\mathcal U}_q({\mathfrak k}^{\mathbb C})$ 
module, \begin{equation}\label{166}{\mathcal M}(V_\Lambda)={\mathcal
A}_q^-\otimes V_\Lambda. \end{equation}

\medskip

We now consider pairings between ${\mathcal M}(V_\Lambda)$ and ${\mathcal
A}_q^+\otimes V^o_{\Lambda'}$.

\medskip

\begin{Def} Let $X=J, K, L$.
$\forall u^+_z\in {\mathcal A}_q^+,u^-_w\in
{\mathcal A}_q^-, v\in V_\Lambda, v'\in V^o_{\Lambda'}$:
\begin{equation}
(u^-_w v,u^+_zv')_X:=(u^-_w ,u^+_z)_X( v',v)_\Lambda.
\end{equation}
Here, $(v',v)_\Lambda$ denotes the natural pairing between a module and 
its dual, and, as usual, the definition is extended by bilinearity to 
the full spaces.\end{Def}
\begin{Rem}There seems to be some bias with this notation when $X=J$. However, 
we use the new form only when considering dual modules. Otherwise, it is the 
form $(\cdot,\cdot)_J$ in Proposition~\ref{p612} that is considered. 
Furthermore, the two forms agree on ${\mathcal A}_q^-\times {\mathcal A}_q^+$. 
Hence we use the symbol $J$ for both forms.
\end{Rem}

\bigskip

By symmetry, the vector space ${\mathcal
A}_q^+V^o_{\Lambda'}$ is also a left module for
${\mathcal U}_q({\mathfrak g})$. {\bf We will, however,
consider another module structure on this space.}

 We extend the notation from Subsection~\ref{5.3} as follows:
\begin{equation}
(\widehat{T}^X(u^+_z v),u^-_wv')_X =(u^+_z v,T(u^-_wv'))_X.
\end{equation}
Recall that $S$ denotes the antipode. We then define
 
\begin{equation}\forall u\in{\mathcal U}_q:\left({\mathcal 
O}_{-\omega(\Lambda)}^X(u)(u^+_z v),u^-_wv'\right)_X=\left(u^+_z 
v,L_\Lambda(S(u))(u^-_wv')\right)_X.\end{equation} In other words,
\begin{equation}{\mathcal 
O}_{\Lambda'}^X(u)=\widehat{{L_\Lambda(S(u))}}^X.\end{equation}

\bigskip

The pairings we consider result in duals of
multiplication operators of the general form :
\begin{equation}\label{gammaX}
{\widehat{M^o}}^X=(H^o)^{\gamma_X}D^o. \end{equation}  Notice that 
$\gamma_J=1$, 
$\gamma_L=0$,
and $\gamma_K=-1$.

\bigskip

\bigskip

\subsection{Relation to holomorphically induced modules}
Let us agree to write \begin{equation}{\mathcal U}_q={\mathcal A}^+_q{\mathcal 
U}_q({\mathfrak k}){\mathcal A}^-_q,\end{equation}and consider functions 
$f:{\mathcal 
U}_q\rightarrow V'_\Lambda$, where $V'_\Lambda$ is a finite dimensional 
${\mathcal U}_q({\mathfrak 
k})$ module. The subspace ${\mathcal H}(V'_\Lambda)$ of such functions which 
furthermore satisfy\begin{eqnarray}
           \forall u\in {\mathcal U}_q, \forall u_{\mathfrak k}\in {\mathcal 
U}_q({\mathfrak 
k}): f(uu_{\mathfrak k})&=&S(u_{\mathfrak k})f(u)\\
\forall u\in {\mathcal U}_q,\forall u_w\in {\mathcal A}^-_q:f(uu_w)&=&0
          \end{eqnarray}
is invariant under left action, and this module is equivalent to our dual 
module. The second condition can be interpreted as saying that the function 
should be annihilated by all anti-holomorphic (quantized) vector fields.

\bigskip

\bigskip

\section{The actions in the module and its dual}

We consider ${\mathcal M}(V_\Lambda)$ and its dual. We only consider $E_\beta$, 
but a similar approach will work for any $E_\gamma$, $\gamma\in\Pi_c$.

\medskip

\subsection{Technical material.} For use in the computation of the dual 
representation, we need to analyze in greater detail the commutator between 
$Z_{11}=E_\beta$ and an element of ${\mathcal U}^+_q({\mathfrak k}_L^{\mathbb 
C})$:

\medskip

We use (cf. (\ref{w-mu}))  \begin{equation}\omega_0 
^L=s_{\mu_{n}}s_{\mu_{n-1}}s_{\mu_{n}}s_{\mu_{n-2}}s_{\mu_{n-1}}s_{\mu_{n
} } \cdot \cdots\cdot s_{\mu_{1}}s_{\mu_{2}}\cdots
s_{\mu_{n}}.\end{equation}

This leads to a PBW basis based on the ordered elements 
\begin{equation}E_{\mu_n}, (T_{\mu_{n}}(E_{\mu_{n-1}})),\dots
(T_{(\omega_0^Ls_{\mu_{n}})}(E_{\mu_{n}})).\end{equation}

Let $u_{L}^+\in {\mathcal U}^+_q({\mathfrak k}_L^{\mathbb C})$. We are 
interested in terms of the form $u_L^+E_\beta$, where we want to move $E_\beta$ 
to the left.

In a straightforward manner, using
Lemma~\ref{156} repeatedly, the right most element in our basis can be seen  to 
be $E_{\mu_{1}}$. Indeed,
the right hand tail of the basis can is 
\begin{equation}T_{\mu_{n}}T_{\mu_{n-1}}\cdots
T_{\mu_{2}}(E_{\mu_{1}}), T_{\mu_{n-1}}\cdots
T_{\mu_{2}}(E_{\mu_{1}}), \dots,
T_{\mu_{2}}(E_{\mu_{1}}), E_{\mu_{1}}.\end{equation}

We begin our computation by observing the equation
\begin{equation}(E_{\mu_{1}})^kE_\beta=[k]T_{\mu_{1}}(E_\beta)E_{\mu_{1}}^{k-1}
+q^{-k}
E_\beta E_ { \mu_ {1}}^k.\end{equation}

Since $T_{\mu_{i}}(E_\beta)=E_\beta$ for $i\geq2$, it
follows that

\begin{equation}T_{\mu_{2}}(E_{\mu_1})E_\beta=T_{\mu_{2}}T_{\mu_{1}}(E_\beta)+q^
{-1}E_\beta
T_{\mu_{2}}(E_{\mu_1}).\end{equation}

Set $\check X_{{i}}=T_{\mu_{i}}T_{\mu_{n-1}}\cdots
T_{\mu_{2}}(E_{\mu_{1}})$ for $i\geq2$, $\check X_1=E_{\mu_1}$, and $\check 
X_0=1$. Using the Serre relations and using that the operators $T_\alpha$ are 
automorphisms, one easily obtains, cf. a similar  computation below, 

\begin{equation}\check X_{{i}}^{a_{i}}E_\beta=[a_{i}]Z_{i+1,1}
\check X_{i}^{a_{i}-1} +q^ { -k } E_\beta
\check X_{{i}}^{a_{i}}. \end{equation}

From the Serre-relations we get (as in (\ref{80})):

\begin{equation}
E_{\mu_k}Z_{a,1}=\left\{\begin{array}{ll}Z_{a+1,1}+qZ_{a,1}E_{\mu_k}
&\textrm{
If }a=k\\q^{-1}Z_{a,1}E_{\mu_k}&\textrm{ If
}a=k+1\\Z_{a,1}E_{\mu_k}&\textrm{ If }a\neq k,
k+1\end{array}\right. .\end{equation}

From this it  follows that \begin{eqnarray}
\check X_{k}Z_{k+1,1}&=&T_{\mu_k}\dots T_{\mu_2}(E_{\mu_1})T_{\mu_k}\dots 
T_{\mu_2}(T_{\mu_1}(E_\beta))\\\nonumber=T_{\mu_k}\dots 
T_{\mu_2}(E_{\mu_1}T_{\mu_1}(E_\beta))&=&qT_{\mu_k}\dots 
T_{\mu_2}(T_{\mu_1}(E_\beta)E_{\mu_1})\\&=&qZ_{k+1,1}\check X_{k}\nonumber
. \end{eqnarray}More generally then,
\begin{equation}
\check X_\ell^kZ_{1,1}=[k]Z_{\ell+1,1}\check X_\ell^{k-1}+q^{-k}Z_{1,1}X_\ell^k.
\end{equation}

Similarly, we get ($i\geq1$)\begin{eqnarray}
\check X_{i+1,1}Z_{i+1,1}&=&(E_{\mu_{i+1}}\check X_{i,1}-q^{-1}\check 
X_{i,1}E_{\mu_{i+1}+})Z_{i+1,1}\\&=&
(q-q^{-1})Z_{i+2,1}\check X_{i,1}+Z_{i+1,1}\check X_{i+1,1}.\end{eqnarray}

There are many more formulas that can be derived in the same manner: 
\begin{equation}\check X_kZ_{k+1}=qZ_{k+1,1}\check 
X_k,\
\check X_kZ_{k+2,1}=Z_{k+2,1}\check X_k,\ \check X_1\check X_2=q^{-1}\check 
X_2\check X_1,\end{equation}
\begin{equation}\check X_2^kZ_{2,1}=(q-q^{-1})Z_{3,1}\check X_2^{k-1}\check 
X_1+Z_{2,1}\check X_2^k,\textrm{ and}\end{equation}
\begin{equation}\check X_2^kZ_{1,1}=Z_{3,1}\check X_2^{k-1}+q^{-k}Z_{1,1}\check 
X_2^k,\end{equation}

More generally, let $i\geq k\geq1$ and $\ell\geq1$. Then
\begin{eqnarray} \check X_{k+\ell,1}Z_{i+1,1}=
(q-q^{-1})Z_{k+\ell+1,1}\check X_{k,1}+Z_{i+1,1}\check X_{k+\ell ,1}.
\end{eqnarray}

It follows that, setting $X_1:=E_{\mu_1}$, \begin{eqnarray}\nonumber
(\check X_k^{a_k}\dots
\check X_i^{a_i}\dots \check 
X_1^{a_1})Z_{1,1}=\sum_{i=1}^kq^{-a_1-\dots-a_{i-1}}
Z_{i+1,1}[a_i](\check X_k^{a_k}\dots \check X_i^{a_i-1}\dots \check X_1^{a_1})
\\ +\ q^{-\sum_ia_i}Z_{1,1}(\check X_k^{a_k}\dots \label{exp}
\check X_i^{a_i}\dots
\check X_1^{a_1})+(q-q^{-1})\cdot(\star\star\star)\end{eqnarray} Above, we 
interpret $q^{-a_1-\dots-a_{i-1}}$ as $1$ when $i=1$.  The expression 
$(\star\star\star)$ in (\ref{exp})  is well behaved under $q\rightarrow
1$ and each monomial term in it is of the same total
degree $1+\sum_{i=1}^ka_i$ with exactly one factor of the form $Z_{r,1}$ for 
some $r=2,\dots,n$. For use in  Corollary~\ref{cor-res} we observe that it thus 
will be unaffected by the change of basis.

We let $\check Y_j$ denote the analogous terms in ${\mathcal U}^+_q({\mathfrak 
k}_R^{\mathbb C})$.
\medskip

\medskip

We have in particular obtained a description of, to what extent $\ad( E_\beta 
K_\beta^{-1})$,  acting on ${\mathcal U}_q({\mathfrak k}^{\mathbb C})$, can 
yield a component in ${\mathcal A}_q^+$. This will be useful later.

\begin{Cor}\label{k-cor}
If $E_\beta Ad(K_\beta^{-1})(\check X_\nu^{\mathbf c}\check Y_\mu^{\mathbf d}) 
- \check X_\nu^{\mathbf c}\check Y_\mu^{\mathbf d}E_\beta \in{\mathcal A}_q^+$ 
then $\check X_\nu^{\mathbf c}\check Y_\mu^{\mathbf d}=\check X_i\check Y_j$ 
for some $i\geq0$ and some $j\geq0$.
\end{Cor}

\medskip

Another easy result, which we will need later, is

\begin{Cor}
\begin{eqnarray}
\forall i=1,\dots, n-1: X_iZ_{11}&=&Z_{i+1,1}+q^{-1}Z_{11}X_i\\
\forall j=1,\dots, n-1: Y_jZ_{11}&=&Z_{1,j+1}+q^{-1}Z_{11}Y_j\\
\forall i,j=1,\dots, n-1: 
X_iY_jZ_{11}&=&Z_{i+1,j+1}
+q^{-1}Z_{i+1,1}Y_j\\&+&q^{-1}Z_{1,j+1}X_i+q^{-2}Z_{11}X_iY_j.
\end{eqnarray}

\end{Cor}

\subsection{$E_\beta$ acting in ${\mathcal U}_q^-$, ${\mathcal 
M}(V_\Lambda)$, and ${\mathcal M}'(V_\Lambda)$.}

 First we introduce PBW bases $X_L^{\mathbf a}\in 
{\mathcal U}_q^-({\mathfrak k}_L^{\mathbb C})$ and $Y_R^{\mathbf b}\in 
{\mathcal U}_q^-({\mathfrak k}_R^{\mathbb C})$ by 
the same recipe as above.  The symbol ${\mathbf a}$, and similarly for the 
others, as usual stands for an $n-1$ touple of non-negative integers. We will 
often 
denote the $n-1$ touple with all zeros as $0$ and the corresponding term 
$X_\nu^{\mathbf a}$  will simply be denoted $X_\nu^{0}$ - which is a 
complicated, but convenient, way of writing $1$.

We maintain the identification $E_\beta=Z_{11}$. To compute the dual action of 
left multiplication by $E_\beta$ in ${\mathcal U}^-_q({\mathfrak g}^{\mathbb 
C})$, we need to compute $\widehat{T}^J$ for $T={{-K^{-1}_\beta LM_{Z_{11}}}}$. 
Here  we use the form $(\cdot,\cdot)_{J}$ from Proposition~\ref{p612} in 
combination with (\ref{116}). First we determine $[E_\beta ,g(w)u^-_{\mathfrak 
k}]=[E_\beta ,g(w)]u^-_{\mathfrak k}$ for $g(w)\in{\mathcal A}_q^-$. It turns 
out to contain all the needed information.

 Set
 \begin{equation}\label{sum-am}
[E_\beta ,g(w)]=\sum_{{\mathbf a},{\mathbf b},k}f^g_{{\mathbf a},{\mathbf 
b},k}X_L^{\mathbf a}Y_R^{\mathbf b}(K_\beta^k)p_{{\mathbf a},{\mathbf b},k} 
 \end{equation}
 where $p_{{\mathbf a},{\mathbf b},k}\in {\mathcal U}^0_q({\mathfrak 
k}_L^{\mathbb C}\oplus {\mathfrak k}_R^{\mathbb C})$, and $f^g_{{\mathbf 
a},{\mathbf b},k} \in{\mathcal A}_q^-$. The summation is over all multiindices 
${\mathbf a},{\mathbf b}$ and all integers $k$ (though it is clear that only 
$k=-1,0,1$ will give rise to something non-zero). It will turn out below that 
if some $p_{{\mathbf a},{\mathbf b},k}\neq0$ then, up to multiplication by a 
non-zero constant,  $p_{{\mathbf a},{\mathbf b},k} =K_\beta^{\ell}$ for some 
integer $\ell=\ell({\mathbf a},{\mathbf b},k)$. We make this choice and thereby 
remove the ambiguity in (\ref{sum-am}).
 Then, noting that $r_\beta$ and $r'_\beta$ map into ${\mathcal U}^-_q$, and 
noting Corollary~(\ref{K}), we have for all $\check p_{{\mathbf c},{\mathbf 
d}}\in {\mathcal U}^0_q({\mathfrak 
k}_L^{\mathbb C}\oplus {\mathfrak k}_R^{\mathbb C})$,  and $\psi_{{\mathbf 
c},{\mathbf 
d}}(z)\in {\mathcal A}^+_q$: \begin{eqnarray}
 ([E_\beta ,g(w)],\psi_{{\mathbf c},{\mathbf d}}(z)\check X_L^{\mathbf 
c}\check Y_R^{\mathbf d}\check p_{{\mathbf c},{\mathbf d}})=\\ 
 \frac1{\gamma}\left(K_\beta r_\beta(g(w))-r'_\beta(g(w))K_\beta^{-1}, 
\psi_{{\mathbf c},{\mathbf d}}(z)\check X_L^{\mathbf c}\check Y_R^{\mathbf 
d}\check p_{{\mathbf c},{\mathbf d}}\right)=\\\label{193}
 \frac1{\gamma(E_\beta,F_\beta)_J}
 \left(g(w),E_\beta Ad(K_\beta^{-1})(\psi_{{\mathbf c},{\mathbf d}}(z)\check 
X_L^{\mathbf c}\check Y_R^{\mathbf d}\right)( K_\beta,
 \check p_{{\mathbf c},{\mathbf d}})\\-\ \frac1{\gamma(E_\beta,F_\beta)_J}
 \left(g(w),\psi_{{\mathbf c},{\mathbf d}}(z)\check X_L^{\mathbf c}\check 
Y_R^{\mathbf d}E_\beta\right)\label{194}
 (K^{-1}_\beta,\check p_{{\mathbf c},{\mathbf d}}).
 \end{eqnarray}

 In (\ref{194}) it follows from Corollary~\ref{k-cor} that for this to be 
non-zero, $X_R^{\mathbf c}\check Y_L^{\mathbf d}=\check X_i\check Y_j$ for 
some $i,j\geq0$. Using the same Corollary, we obtain
 \begin{eqnarray}&([E_\beta ,g(w)],\psi_{{\mathbf c},{\mathbf d}}(z)\check 
X_\nu^{\mathbf c}\check Y_\mu^{\mathbf d}\check p_{{\mathbf c},{\mathbf 
d}})=&\\&-
\delta_{({\mathbf c},{\mathbf d}), (0,0)} \left(g(w),E_\beta 
Ad(K_\beta^{-1})(\psi_{0,0})(z)\right)+&\\&
  \delta_{({\mathbf c},{\mathbf d}), 
(0,0)}\left(g(w),(\psi_{0,0})(z)E_\beta\right)
+\sum_{i+j>0}\delta_{({\mathbf c},{\mathbf d}), 
(e_i,e_j)}\left(g(w),\psi_{ij}Z_{i+1,j+1}\right)\nonumber=&\\&-
 \delta_{({\mathbf c},{\mathbf d}), (0,0)}\left(g(w),\left({}_{Z_{11}}M\right)
Ad(K_\beta^{-1})(\psi_{0,0})(z)\right)(K_\beta,\check p_{0,0})+\label{198}&\\&
  \delta_{({\mathbf c},{\mathbf d}), (0,0)}\left(g(w),(M_{Z_{11}}( 
\psi_{0,0}))\right)(K^{-1}_\beta,\check p_{0,0})\label{199}
+&\\&\sum_{i+j>0}\delta_{({\mathbf c},{\mathbf d}), 
(e_i,e_j)}\left(g(w),M_{Z_{i+1,j+1}} \psi_{ij}\right)(K^{-1}_\beta,\check 
p_{i,j}).&\label{200}
 \end{eqnarray}
 
 On the other hand,\begin{eqnarray}
 ([E_\beta ,g(w)],\psi_{{\mathbf c},{\mathbf d}}(z)\check X_\nu^{\mathbf 
c}\check Y_\mu^{\mathbf d}\check p_{{\mathbf c},{\mathbf 
d}})=\\\label{202}\left(\sum_{{\mathbf a},{\mathbf b},k}f^g_{{\mathbf 
a},{\mathbf b},k}X_\nu^{\mathbf a}Y_\mu^{\mathbf b}(K_\beta^k)p_{{\mathbf 
a},{\mathbf b},k} ,\psi_{{\mathbf c},{\mathbf d}}(z)\check X_\nu^{\mathbf 
c}\check Y_\mu^{\mathbf d}\check p_{{\mathbf c},{\mathbf d}}\right),
 \end{eqnarray}
 where clearly $({\mathbf a},{\mathbf b})=({\mathbf c},{\mathbf d})$. 
 
 \smallskip

 In (\ref{198}-\ref{200}) we may introduce  the duals of the left and 
right multiplication operators, whereby one can compare directly to 
(\ref{202}). These duals are here denoted   $O_1,O_2,O_{ij}: {\mathcal 
A}_q^-\rightarrow 
{\mathcal A}_q^-$. 

\smallskip

Using the non-degeneracy of the form, we then reach:

\begin{Thm}\label{main1}
\begin{eqnarray}\nonumber[E_\beta,g(w)]&=&O_1(g)K_\beta + O_2(g)K_\beta^{-1}\\
&+&\sum_{i=1}^{n-1}O_{i1}(g)\frac{X_i}{(X_i,\check X_i)_J}K_\beta^{-1}
+\sum_{j=1}^{n-1}O_{1j}(g)\frac{Y_j}{(Y_j,\check 
Y_j)_J}K_\beta^{-1}\\&+&\sum_{i\cdot j>0}O_{ij}(g)\frac{X_iY_j}{(X_i,\check 
X_i)_J(Y_j,\check Y_j)_J}K_\beta^{-1}, \nonumber\textrm{ where}\end{eqnarray}
\begin{eqnarray}\forall g, \forall \psi: 
(O_1(g),\psi)_J&=&-(g,Z_{11}Ad(K_\beta^{-1})(\psi))_J,
\\\forall g, \forall \psi: (O_2(g),\psi)_J&=&(g,\psi_{00}Z_{11})_J,\\\forall g,
\forall i,j; \forall \psi: (O_{ij}(g),\psi)_J&=&(g,\psi 
Z_{i+1,j+1})_J.\end{eqnarray}
\end{Thm}

\medskip

Recall that  that $\frac{1}{(X_i,\check X_i)_J}=\frac{1}{(Y_j,\check 
Y_j)_J}=-\gamma$.

\medskip

\begin{Cor}\label{actJ}
\begin{eqnarray}{\mathcal O}^J_{\Lambda'}(E_\beta)=
\left({}_{Z_{11}}M\right)\otimes 1 - M_{Z_{11}}Ad(K_\beta)\otimes K_\beta^2 
\\\nonumber
+\sum_{i>0}\gamma q^{-1}M_{Z_{i+1,1}}Ad(K_\beta)\otimes 
X_i^TK_\beta^2+  \sum_{j>0}\gamma 
q^{-1}M_{Z_{j+1,1}}Ad(K_\beta)\otimes Y_j^TK_\beta^2\\\nonumber 
-\sum_{i\cdot j>0}\gamma^2q^{-2}M_{Z_{i+1,j+1}}Ad(K_\beta)
\otimes X_i^TY_j^TK_\beta^2,
\end{eqnarray}
where the operators $X_i^T$ and $Y_j^T$ are the dual operators to $X_i,Y_j$ as 
acting in $V_\Lambda$.
\end{Cor}

\bigskip

We now address the issues of rescaling to avoid the singularity of the form at 
$q=1$.

\smallskip

We let ${\mathbb B}(W^{\mathbf b})=(-\gamma)^{\mid{\mathbf
b}\mid}W^{\mathbf b}$, but do not rescale in $V_\Lambda$ or its dual.

\begin{Def}
\begin{equation}\forall u\in{\mathcal U}_q:
{\mathcal O}^{{\mathbb B},J}_{\Lambda'}(u)={\mathbb B}^{-1}{\mathcal 
O}_{\Lambda'}^J(u){\mathbb
B}.
\end{equation}
\end{Def}

\begin{Cor}
\begin{equation}{\mathbb B}^{-1}{\mathcal O}_{\Lambda'}^J(F_\beta){\mathbb
B}=-\left(Ad(K_\beta^{-1}\right)\circ
r'_\beta.\end{equation}
	\end{Cor}This has a well-defined limit at $q=1$.

\medskip

\begin{Cor}\label{cor-res}
\begin{eqnarray}
{\mathcal O}^{{\mathbb B},J}_{\Lambda'}(E_\beta)=-\frac1{\gamma} \left(
\left({}_{Z_{11}}M\right)\otimes 1 -M_{Z_{11}}Ad(K_\beta)\otimes 
K_\beta^2\right)\\ 
\quad {-} \sum_{i>0} q^{-1}M_{Z_{i1}}Ad(K_\beta)\otimes X_i^TK_\beta^2 
{-} \sum_{j>0} q^{-1}M_{Z_{j1}}Ad(K_\beta)\otimes Y_j^TK_\beta^2\\ + 
\sum_{i\cdot j>0}\gamma q^{-2}M_{Z_{ij}}Ad(K_\beta)
\otimes X_i^TY_j^TK_\beta^2.
\end{eqnarray}
\end{Cor}

\medskip

\subsection{Comparison to the classical result}
The following is well known, see e.g. (\cite{jv}): If 
$X=\left(\begin{array}{cc}X_1&X_2\\X_3&X_4\end{array}\right)\in 
su(n,n)^{\mathbb C}$, the
infinitesimal action is given as
\begin{eqnarray}\label{a-act}
dU_\pi(X)f(z)&=&d\pi\left(\begin{array}{cc}X_1-zX_3&0\\0&X_3z+X_4\end{array}
\right)f(z)\\&-&
(\delta(X_1z+X_2-zX_3z-zX_4)f)(z).\nonumber
  \end{eqnarray}
Here, $d\pi$ is the dual of the finite-dimensional representation of ${\mathcal 
U}_q({\mathfrak k}^{\mathbb C})$ of highest weight 
$\Lambda=(\Lambda_L,\Lambda_R,\lambda)$. The version of $su(n,n)^{\mathbb C}$ 
we use is the one based on the Hermitian form 
$H_\beta=\left(\begin{array}{cc}I_n&0\\0&-I_n\end{array}\right)$. Finally, 
$\delta(Y)$ denotes the directional derivative 
$\sum_{ij}Y_{ij}\frac{\partial}{\partial Z_{ij}}$.

\bigskip

We now cite a special case of a result which will be proved in \S\ref{multiop}:

\begin{Lem}
\begin{eqnarray}&M_{Z_\beta}=(H^0_{n,1}\dots
H^0_{2,1}H^0_{1,2}\cdots H^0_{1,n})M_{11}^0+\\\nonumber
&\gamma\sum_{c=2,d=2}^{n,n}P_{c,d}(H^0)\left(D^o_{c,d}M^o_{c,1}M^o_{1,d} 
\right)+
\gamma^2T_2 \end{eqnarray}
for monomials $P_{c,d}(H^0)$ and some term $T_2$  proportional to at least 
$\gamma^2$.
\end{Lem}

It follows easily that 
\begin{eqnarray}
{\mathcal O}^{{\mathbb B},J}_{\Lambda'}(E_\beta)(Z^{\mathbf 
a})=-\frac1{\gamma}\left(1-q^{2(a_{11}+a_{21}+\cdots
+a_{n1}+a_{12}+\cdots
+a_{1n}-2\lambda)}\right)Z_{11}Z^{{\mathbf a}}\\+\sum_{i=2,j=2}^{n,n}q^{
\alpha( { \mathbf a},i,j,\lambda)} D^o_{i,j}M^o_{i,1}M^o_{1,j}(Z^{\mathbf a})
 {-} \sum_{i>0} q^{-1}M_{Z_{i1}}Ad(K_\beta)\otimes X_i^TK_\beta^2\\ 
{-} \sum_{j>0} q^{-1}M_{Z_{j1}}Ad(K_\beta)\otimes Y_j^TK_\beta^2+
\gamma T_2.\nonumber
\end{eqnarray}

If we let $Z_1^{\underline{a}}$ denote the limit as
$q\rightarrow 1$ then
\begin{eqnarray}\nonumber\lim_{q\rightarrow1}\frac1{\gamma}\left(1-q^{2(a_{11}
+a_{21}+\cdots
+a_{n1}+a_{12}+\cdots
+a_{1n}-2\lambda)}\right)Z_{11}Z^{{\mathbf a}}Z^{\underline{a}}=\\(-Z^2_{11}
\frac{\partial}{\partial
Z_{11}}-Z_{11}Z_{21}\frac{\partial}{\partial
Z_{21}}-\cdots-Z_{11}Z_{n1}\frac{\partial}{\partial
Z_{n1}}\\-Z_{11}Z_{12}\frac{\partial}{\partial
Z_{12}}-\cdots -Z_{11}Z_{1n}\frac{\partial}{\partial
Z_{1n}}+Z_{11}\lambda\cdot I)Z_1^{{\mathbf a}}.\end{eqnarray}

\medskip

\begin{Cor}\label{class-cor}
\begin{eqnarray} lim_{q\rightarrow 1}{\mathcal O}^{{\mathbb 
B},J}_{\Lambda'}(E_\beta)=\sum_{i,j=1}^nZ_{i1}Z_{ij} 
\frac{\partial}{\partial Z_{ij}}- Z_{11}\lambda\cdot I\\
{-} \sum_{i>0} Z_{i1}\otimes X_i^T {-} \sum_{j>0} 
Z_{j1}\otimes Y_j^T.
\end{eqnarray}
\end{Cor}

\medskip

\begin{Obs}
If we agree to write our matrix Z representing the elements $Z_{ij}$ as 
$\sum_{ij}Z_{ij}E_{n+1-i,j}$, the matrix $W$ representing the $W_{ij}$ as 
$\sum_{ij}W_{ij}E_{i,n+1-j}$ and $X_3=E_{1,n}$, then we get that the formula in 
Corollary~\ref{class-cor} is the same as (\ref{a-act}): Recall that we work 
with the dual representation on the 
${\mathfrak k}$ level. Specifically, $X_i=E_{n, n-i}$ and $Y_j=-E_{j+1,1}$ are 
exactly the correct expressions, bearing in mind the way the elements $W_{ij}$ 
are defined, cf. 
(\ref{9}) and (\ref{26}). Also notice that 
\begin{equation}-Z_{11}\lambda=(Z_{11}(-\Lambda(\beta)),\end{equation}as must 
be the case in the dual module.
\end{Obs}

\medskip

\bigskip \section{Differential operators}

\bigskip

\subsection{Multiplication operators and their duals}We now introduce the 
fundamental multiplication operators. Together with their duals they form the 
foundation of any reasonable algebra of differential operators. Here we will 
be interested in examining how the different pairings may lead to different 
algebras. It should be noted that we use the ``bare duality'' to 
define differential operators. By this we mean that we do not use the antipode 
from any Hopf algebra that may be otherwise naturally affiliated with the 
situation.

\smallskip

\begin{Def}We define linear operators \begin{eqnarray} 
\left({}_{{W_{ij}}}M\right),
M_{{W_{ij}}}:{\mathcal U}_q^-&=&{\mathcal
A}_q^-{\mathcal U}_q^-({{\mathfrak k}^\mathbb C})\rightarrow
{\mathcal U}_q^- \textrm{ by}\\
{\mathcal A}_q^-{\mathcal
U}_q^-({\mathfrak k}^{\mathbb C})\ni u_Wu_{{\mathfrak k}^\mathbb 
C}^-&\stackrel{M_{W_{i,j}}}{\rightarrow}&u_WW_{i,j} u_{{\mathfrak
k}^\mathbb C}^-\\ {\mathcal A}_q^-{\mathcal
U}_q^-({\mathfrak k}^{\mathbb C})\ni u_Wu_{{\mathfrak
k}^\mathbb C}^-&\stackrel{\left({}_{{W_{ij}}}M\right)}{\rightarrow}&W_{i,j}u_W 
u_{{\mathfrak 
k}^\mathbb C}^-.\end{eqnarray}
\label{multiop}
Likewise, we define linear operators $\left({}_{{Z_{ij}}}M\right),
M_{{Z_{ij}}}:{\mathcal U}_q^+={\mathcal
A}_q^+{\mathcal U}_q^+({\mathfrak k}^\mathbb C)\rightarrow
{\mathcal U}_q^+$ by   \begin{eqnarray}
{\mathcal A}_q^+{\mathcal
U}_q^+({\mathfrak k}^{\mathbb C})\ni u_Zu_{{\mathfrak k}^\mathbb 
C}^+&\stackrel{M_{Z_{i,j}}}{\rightarrow}&u_ZZ_{i,j} u_{{\mathfrak
k}^\mathbb C}^+,\\ {\mathcal A}_q^+{\mathcal
U}_q^+({\mathfrak k}^{\mathbb C})\ni u_Zu_{{\mathfrak
k}^\mathbb C}^+&\stackrel{\left({}_{{Z_{ij}}}M\right)}{\rightarrow}&Z_{i,j}u_Z 
u_{{\mathfrak 
k}^\mathbb C}^+.\end{eqnarray}   Recall
 (\ref{41} - \ref{43}).

We define the linear operator
$\widehat{M_{W_{i,j}}}^X$ acting on ${\mathcal
A}_q^+{\mathcal U}_q^+({{\mathfrak k}^\mathbb C})$ by
\begin{eqnarray}( u_Wu_{\mathfrak
k}^-,\widehat{M_{W_{i,j}}}^X(u_Z u_{\mathfrak
k}^+))_X&=&({M_{W_{i,j}}} (u_Wu_{{\mathfrak k}^\mathbb C}^-),u_Z
u_{\mathfrak k}^+)_X\\ &=&((u_WW_{ij}u_{\mathfrak k}^-),u_Z
u_{{\mathfrak k}^\mathbb C}^+)_X. \end{eqnarray} The linear operator
$\widehat{\left({}_{{W_{ij}}}M\right)}^X$ is defined analogously.
\end{Def}
 
 \smallskip
 
 We also use the notation
 
\begin{equation} {\frac{\partial}{\partial
Z_{ij}}}_{X}=\widehat{M_{W_{i,j}}}^X\textrm{ and }
{\vphantom{\int}}_X {\frac{\partial}{\partial
Z_{ij}}}=\widehat{{}_{{W_{ij}}}M}^X. 
 \end{equation}
 
 Similarly, of course, for all the others.
 
From the behavior of the bilinear form, it follows that
\begin{eqnarray}{\frac{\partial}{\partial
Z_{ij}}}_{X}\left(Z^{\mathbf a}u_{\mathfrak
k}^+\right)&=&\left({\frac{\partial}{\partial
Z_{ij}}}_{X}Z^{\mathbf a}\right)u_{\mathfrak k}^+\\
{\vphantom{\int}}_X {\frac{\partial}{\partial Z_{ij}}}\left(Z^{\mathbf
a}u_{\mathfrak k}^+\right)&=&\left({\vphantom{\int}}_X {\frac{\partial}{\partial
Z_{ij}}}Z^{\mathbf a}\right)u_{\mathfrak
k}^+.\end{eqnarray}
 
 \medskip

\subsection{Four algebras of differential operators}
On the way to defining differential operators on
${\mathcal A}_q^{+}$ we make some preliminary
definitions:
   
\begin{Def}Set
\begin{equation}{\mathcal D}(\partial _X,{}_ZM)={\mathbb
C}\left({(\frac{\partial}{\partial
Z_{ij}})}_{X},{}_{Z_{k,\ell}}M\mid
i,j,k,\ell\in\{1,2,\dots,n\}\right);\end{equation}
the algebra  of
quantized right-left differential operators.

Similarly, \begin{eqnarray} {\mathcal
D}(\partial_X,M_Z)&:=&{\mathbb C}\left({(\frac{\partial}{\partial
Z_{ij}})}_{X},M_{Z_{k,\ell}}\mid
i,j,k,\ell\in\{1,2,\dots,n\}\right)\\ {\mathcal D}(
{}_X\partial,{}_ZM)&:=&{\mathbb 
C}\left({\vphantom\int}_X{(\frac{\partial}{\partial
Z_{ij}})},{}_{Z_{k,\ell}}M\mid
i,j,k,\ell\in\{1,2,\dots,n\}\right)\\ {\mathcal
D}({}_X\partial,M_Z)&:=&{\mathbb 
C}\left({\vphantom\int}_X{(\frac{\partial}{\partial
Z_{ij}})},M_{Z_{k,\ell}}\mid
i,j,k,\ell\in\{1,2,\dots,n\}\right) \end{eqnarray} 
with analogous names.\end{Def}

\bigskip

\subsection{Explicit formulas for the left multiplication operators} We will 
now  compute $\left({}_{W_{i,j}}M\right)$ and
$M_{W_{i,j}}$ explicitly. We do so by using (\ref{a} - \ref{cross})
repeatedly. Notice that ${\mathcal A}_q^+$ has the same
relations as  ${\mathcal A}_q^-$, so that we by the same
computations also compute $\left({}_{Z_{i,j}}M\right)$ and
$M_{Z_{i,j}}$.

\medskip

Specifically, we wish to expand
$W_{ij}W_{11}^{a_{11}}\dots W_{xy}^{a_{xy}}\dots
W_{ij}^{a_{ij}}$ in our basis (\ref{20}) (cf.
Proposition~\ref{2.5}) by using these formulas. Unless
$i=1$ or $j=1$ - in which case we can regroup using only
the equations (\ref{a} - \ref{c}) - we have to use
(\ref{cross}), and hence (\ref{2cross}), already at the position $(1,1)$. Let
us be more general and say that we have reached a
position $(x,y)$ where we use (\ref{2cross}) to make the
replacement $W_{ij}W_{xy}^{a_{xy}}\rightarrow
W_{ij}W_{xy}^{a_{xy}}W_{ij}+(q-q^{-1})
q^{a_{xy}-1}[a_{xy}]_qW_{xy}^{a_{xy}-1}W_{xj}W_{iy}$.
Let us focus on the second term: $W_{xj}$ is already in
its right row and can be placed in its correct position
using (\ref{a}). $W_{iy}$ however, may be in a wrong row and, if so, 
to bring it into its correct position we will have to use
(\ref{cross}) again. This means that we have to keep
track of the number of times we use either of the
equations (\ref{a}, \ref{b}), and (\ref{cross}). We omit the details of this 
cumbersome bookkeeping. The
result we obtain is:
 
\begin{Lem}\label{tedi} Let $i,j$ be given. Let
$r\in{\mathbb N}$, and let \begin{equation}(\underline a,\underline
b)=((a_1,a_2,\dots, a_r),(b_1,b_2,\dots,b_r))\in{\mathbb
N}^r\times {\mathbb N}^r.\end{equation} We say that $(\underline
a,\underline b)$ is a NW-partition of $(i,j)$ if
\begin{eqnarray} a_0=1\leq a_1<a_2<\dots
<a_r<a_{r+1}:=i\textrm{ and }\\b_0=1\leq b_1<b_2<\dots
<b_r<b_{r+1}:=j,\end{eqnarray} and we let ${\mathcal
P}^{NW}_r$ denote the set of all such. In the following formulas we use the 
convention that $(\prod_{j\in\emptyset}
F_j)=1$. The following formula holds:
\begin{eqnarray}\label{l-w}&\left({}_{W_{i,j}}M\right)=
M^o_{i,j}(H^o_{i,1}\dots H^o_{i,j-1})(H^o_{i-1,j}\dots
H^o_{1,j})+\\& \gamma^{r}\cdot\sum_{r\in{\mathbb N}}\sum_{(\underline
a,\underline b)\in {\mathcal P}^{NW}_r}
\prod_{k=1}^{r+1}(\prod_{y\in ]b_{r-k+1}, 
b_{r-k+2}[}H^o_{a_k,y})(\prod_{x\in]a_{k-1},a_k[}H^o_{x,b_{r-k+2}})
\cdot\nonumber\\
&M^o_{a_{r+1},b_1}\cdot
\prod_{k=1}^{r}(H^o_{a_k,b_{r-k+1}}D_{a_k,b_{r-k+1}}^o
M^o_{a_{k},b_{r-k+2}})\nonumber .\end{eqnarray} \end{Lem}

We also set

\begin{equation}\label{mdm}B^{NW}_r(\underline a,\underline
b)=M^o_{a_{r{+1}},b_1}\cdot
\prod_{k=1}^{r}(D_{a_k,b_{r-k+1}}^o
M^o_{a_{k},b_{r-k+2}}).\end{equation}
                        When we analyze further on
(\ref{l-w}), the terms in (\ref{mdm}) are important, while the
precise form of the monomials in the elements $H^o_{st}$
are of no importance.

\medskip

The following observation is very useful in later computations:
\begin{Obs}\label{obs1}
If $k<r+1$ and $M^o_{a_k,y_1}$ is a factor in (\ref{mdm}), then so is 
$D^o_{a_k,y_2}$ for precisely one  $y_2$, and $y_2<y_1$. $M^o_{a_{r+1},b_1}$ is 
a factor, but no factor $D^o_{a_{r+1},y}$ occurs in (\ref{mdm}). Similarly, If 
$k<r+1$ and $M^o_{x_1,b_k}$ is a factor in (\ref{mdm}), then so is 
$D^o_{x_2,b_k}$ for precisely one  $x_2$ and $x_2<x_1$. $M^o_{a_1,b_{r+1}}$ is 
a 
factor, but no factor $D^o_{x,b_{r+1}}$ occurs.
\end{Obs}

\medskip

\begin{Rem}
 $\left({}_{Z_{i,j}}M\right)$ is given by exactly the same formula as $\left({}_{W_{i,j}}M\right)$. The operators $M^o,H^o$, and $D^o$ must just be interpreted as operators in ${\mathcal A}_q^+$.
\end{Rem}

\medskip

\subsection{The right multiplication operators}

Here we get a similar result. In fact, using the anti-automorphism 
$\phi:Z_{ij}\rightarrow Z_{n+1-i,n+1-j}$ one gets the expression for 
$M_{W_{ij}}$ from that of $\left({}_{W_{n+1-i,n+1-j}}M\right)$ by replacing all 
factors 
$X_{a,b}$, $X=D^o,M^o,H^o$ by $X_{n+1-a,n+1-a}$ in the latter. We choose to 
state 
it with fewer details since it 
is only the explicit form of the term involving $M^o_{ij}$, together with 
Observation~\ref{obs2} below, that is important.
\begin{Lem}\label{tedi} Let $i,j$ be given. Let
$r\in{\mathbb N}$, and let \begin{equation}(\underline a,\underline
b)=((a_1,a_2,\dots, a_r),(b_1,b_2,\dots,b_r))\in{\mathbb
N}^r\times {\mathbb N}^r.\end{equation} We say that $(\underline
a,\underline b)$ is an SE-partition of $(i,j)$ if
\begin{equation}i=a_0<a_1<\dots<a_r\leq a_{r+1}=n\textrm{ and
}j=b_0<b_1<\dots<b_r\leq b_{r+1}=n\}. \end{equation} We denote
by ${\mathcal P}^{SE}_r$ the set of all such.

For $(\underline{a},\underline{b})\in {\mathcal
P}^{SE}_r$ set
\begin{eqnarray}\label{rw1}B_RM^{SE}(\underline{a},\underline{b})&=&M^o_{a_r,b_0
}
\prod_{x=0}^{r-1}(D^o_{a_{r-x},b_{x+1}}M^o_{a_{r-x-1},b_{x+1}})
\textrm{ and 
}\\C(\underline{a},\underline{b})&=&D^o_{a_r,b_0}\prod_{x=0}^{r-1}(M^o_{a_{r-x},
b_{
x+1}}D^o_{a_{r-x-1},b_{x+1}}) .\end{eqnarray}

Then \begin{eqnarray}\label{rw2}M_{W_{ij}}=\\
M^o_{i,j}(H^o_{i,j+1}\dots H^o_{i,n})(H^o_{i+1,j}\dots
H^o_{n,j})\\+\sum_{r\in{\mathbb
N}}\sum_{(\underline{a},\underline{b})\in {\mathcal
P}^{SE}_r}B_RM^{SE}(\underline{a},\underline{b})H^o_r(\underline{a},\underline{b
}
),\end{eqnarray} where each
$H^o_r(\underline{a},\underline{b})$ is a Laurent monomial
in some of the elements $H^o_{s,t}$.

Furthermore\begin{eqnarray}&\widehat{M_{W_{ij}}}^X=\\&(H^o_{ij})^{\gamma_X}D^o_{i,j}(H^o_{i
, j+1
}\dots
H^o_{i,n})(H^o_{i+1,j}\dots
H^o_{n,j})\sum_{(\underline{a},\underline{b})\in {\mathcal
P}^{SE}_r}C(\underline{a},\underline{b})\widetilde{H^o_r}^X(\underline{a},
\underline{b})\nonumber,\end{eqnarray}where
$\widetilde{H^o_r}^X(\underline{a},\underline{b})=q^{\alpha_X}H^o_r(\underline{a
},
\underline{b})$ and $\alpha_X$ depends on the case (as
well as on (\underline{a},\underline{b})). In the
following analysis, the exact value is of no importance.
\end{Lem}

\medskip

\begin{Obs}\label{obs2}Similarly to Observation~\ref{obs1}, each row either 
has none or has exactly one $D$ and
one $M$ except row $a_{0}=i$. Each column has either none or exactly
one $D$ and one $M$ except column $b_0=j$.\end{Obs}

\bigskip

\subsection{The various algebras of differential operators}

We now combine left or right multiplication and left or right differential 
operators into bigger algebras of differential operators. To begin with, to 
understand their individual significance, we combine one kind of multiplication 
operator with one kind of differential operator. Before embarking on this, we 
need a special tool:

\subsubsection{Eigenspace decompositions in subalgebras}
\begin{Def}\begin{equation} {\mathcal L}:={\mathbb
C}[(H^o_{ij})^{\pm1}\mid i,j=1,\dots,n]\end{equation} denotes the
Laurent polynomial algebra generated by the commuting
invertible elements $H^o_{ij}; i,j=1,\dots,n$. Furthermore, we set
\begin{equation}{\mathcal L}{\mathcal M}=\{\phi\in{\mathcal
L}\mid\phi\textrm{ is a monomial}\}.\end{equation} \end{Def}

\medskip

\begin{Prop}\label{goodprop}Let $=O_0, O_1,\dots,O_N\in
{\mathcal W}eyl_q(n,n)$ and let
$\Phi=\{\phi_1,\dots,\phi_N\}\subset {\mathcal
L}{\mathcal M}$. Assume that \begin{eqnarray}\forall i=1,\dots,N,
\forall j=0,1,\dots,N: \phi_iO_j=q^{k_{ij}}O_j\phi_i\\\Leftrightarrow 
Ad\phi(O_j)=q^{k_{ij}}O_j \textrm{ for some } k_{ij}\in{\mathbb 
Z}.\end{eqnarray}

Suppose that $\Phi$ {\it distinguishes elements} in the following sense: 
$k_{ij}\neq k_{ii}$ whenever
$j\notin\{1,\dots,i\}$.

Let ${\mathcal A}_0$ be a subspace of ${\mathcal
W}eyl_q(n,n)$ which is invariant under right and left multiplication by 
$\phi_x;x=0,1,\dots,N$. Suppose that
$O=\sum_{j=0}^NO_j$ belongs to ${\mathcal A}_0$, then 
\begin{equation}\label{produ}O_0 \left(\prod_{i=1}^N\phi_i\right)\in 
{\mathcal
A}_0.\end{equation} \end{Prop}

\proof We have by assumption that \begin{equation}\phi_1
O-q^{k_{11}}O\phi_1=\sum_{
j\in\{0,2,\dots,N\}}(q^{k_{1j}}-q^{k_{11}})O_j\phi_1\in{\mathcal
A}_0.\end{equation} Up to irrelevant complex factors, we have thus
removed the summand $O_1$. By the further assumptions we
can continue to remove summands until only $O_0$
remains. \qed

\medskip

\begin{Rem}We will use Proposition~\ref{goodprop}
repeatedly in the
sequel. Here, ${\mathcal  A}_0$ will even be a subalgebra, but we are not 
assuming that the inverses $\phi_x^{-1},x=1,\dots,N$ stabilize ${\mathcal  
A}_0$, so we have to keep the factor $\prod_{i=1}^N\phi_i$ in (\ref{produ}).
\end{Rem}

\bigskip

Now we are ready for the algebras:

\subsubsection{Right - Left}\label{rl}$\quad$\newline First we consider ${\mathcal A}_0={\mathcal
D}(\partial_X,{}_ZM)$. Here we have all operators
${}_{Z_{ij}}M$ and $\widehat{M_{W_{ij}}}^X$.

\medskip

First observe the following simple fact which follows
from formulas (\ref{f1}-\ref{f2}):
 
\begin{Lem}If ${\mathcal D}(\partial_X,{}_ZM)$ contains
elements $\phi_{ij}M^o_{i,j}$ and $\psi_{ij}D^o_{i,j}$,
with $\phi_{ij},\psi_{ij}$ elements of ${\mathcal
L}{\mathcal M}$, then \begin{equation}\label{2}
\phi_{ij}\psi_{ij}(H^o_{i,j})^{\pm1}\in {\mathcal
D}(\partial_X,{}_ZM). \end{equation} \label{lem-pm} \end{Lem}
 
 \smallskip
 
The signs in (\ref{2}) is the source of an ambiguity which
we try to control somewhat by the specific choices below
(\ref{bias}).

\medskip

We now begin to prove that for all $i,j\in\{1,2,\dots,n\}$,
${\mathcal D}(\partial_X,{}_ZM)$ indeed does contain elements
$\phi_{ij}M^o_{i,j}$ and $\psi_{ij}D^o_{i,j}$ and
$\phi_{ij},\psi_{ij}$ may be explicitly given in ${\mathcal L}{\mathcal M}$. 

\medskip

First observe, for use here and later, that by definition, 
\begin{equation}\nonumber\left({}_{Z_{1,j}}M\right) =H^o_{1,1}\dots 
H^o_{1,j-1}M^o_{1,j} \in
{\mathcal D}(\partial_X,{}_ZM)\end{equation} and $\widehat{RM^-_{W_{i,n}}}^X
\in {\mathcal D}(\partial_X, {}_ZM)$ where
$M_{W_{i,n}}=H^o_{nn}\dots H^o_{i+1,n}M^o_{i,n}$.  At the position $(1,n)$ 
we thus
have \begin{equation}\label{first}H^o_{1,1}\cdots
H^o_{1,n-1}M^o_{1,n}, \ H^o_{nn}\cdots
H^o_{2,n}(H^o_{1,n})^{\gamma_X}D^o_{1,n} \in {\mathcal 
D}(\partial_X,{}_ZM),\end{equation} and then, immediately from
Lemma~\ref{lem-pm},
\begin{equation}\label{second}H^o_{1,1}\cdots
H^o_{1,n-1}\cdot H^o_{nn}\cdots H^o_{2,n}(H^o_{1n})^{\gamma_X\pm1}
\in {\mathcal D}(\partial_X,{}_ZM).\end{equation}

\medskip There are 2 covariant first order elements in
our quadratic algebra, namely $Z_{1,n}$ and $Z_{n,1}$.
We have chosen to construct elements working out from
 position $(1,n)$, but of course one may as well use
position $(n,1)$, or a mix of the two.

We use an ordering \begin{equation}(1,n)<(2,n)<\dots
<(n,n)<(1,n-1)<(2,n-1)<\cdots \ \cdots<(n,1),\end{equation} but we
still need some definitions before we get to a precise
statement:

\begin{Def}Introduce the following elements in ${\mathcal L}{\mathcal M}$: 
\begin{eqnarray}V^{NW}_{ij}&=&(\prod_{k<
i}H^o_{kj}) (\prod_{s<j}H^o_{is}) ,\\
V^{SE}_{ij}&=&(\prod_{k>i}H^o_{kj}) (\prod_{s>j}H^o_{is})
,\textrm{ and}\\\label{bias} V_{ij}&=&
V^{NW}_{ij}V^{SE}_{ij}(H^o_{ij})^{\alpha_X},\ \alpha_X=
\left\{\begin{array}{l}0\textrm{ in case }J\\-2\textrm{
in case }K\\ -1\textrm{ in case }L\end{array}\right. .
\end{eqnarray} Inductively set \begin{equation}{\mathbb
V}_{in}:=\left(\prod_{s<i}{\mathbb
V}_{sn}\right){V}_{in}=\left(\prod_{s<i}{V}^{2^{i-1-s}}_{sn}\right){V}_{in}.\end
{equation} More generally, set 
\begin{equation}{\mathbb
V}_{ij}=V_{ij}\left(\prod_{a=1}^{i-1}{\mathbb
V}_{aj}\right)\left(\prod_{b=j+1}^{n}{\mathbb
V}_{i,b}\right). \end{equation} In particular,
\begin{equation}{\mathbb
V}_{1j}:=\left(\prod_{s=1}^{(n-j)}
V_{1,j+s}^{2^{s-1}}\right) V_{1j}.\end{equation}
\end{Def}

\begin{Rem}  Below, we use heavily Observation~\ref{obs1}
and Observation~\ref{obs2}. Our main tool is Proposition~\ref{goodprop} and 
here we need to construct a family $\Phi$  in ${\mathcal L}{\mathcal M}$ that  
distinguishes elements, e.g. 
$D^o_{a,j}M^o_{ij}M^o_{in}$, $a=1,\dots,i-1$ from
$M^o_{in}$. For this purpose we need to assume $\alpha_X\neq 2$. With 
$\gamma_X$ as in (\ref{gammaX}), we have
$\alpha_X=\gamma_X\pm1$, so this  only poses restrictions in case $X=J$ where 
we must pick $\alpha_X=0$. In the other cases there is the previously mentioned 
ambiguity at each position $i,n$.
This means that the elements constructed in the following lemmas are not 
unique, 
but they suffice for our purposes. Furthermore, it does not seem to give 
simplifications in the end results by allowing more general vales. The above 
choices are then made for the sake of specificity.\end{Rem}

 \medskip

\begin{Lem} It holds that \begin{equation}\left(\prod_{s<i}{\mathbb
V}_{sn}\right)V^{NW}_{i,n}M^o_{in} \in {\mathcal 
D}(\partial_X,{}_ZM)\textrm{  and }\widehat{M_{W_{in}}}^X \in
{\mathcal D}(\partial_X,{}_ZM),\end{equation} where
$M_{W_{in}}=V_{in}^{SE}M^o_{in}$. Hence, in particular, the elements ${\mathbb 
V}_{i,n}$ 
belong to ${\mathcal D}(\partial_X,{}_ZM)$.  \end{Lem}

\proof This is proved by induction. The case $i=1$ is
just (\ref{first}) and (\ref{second}) with special
choices.

Now look at a position $(i,n)$ with $i\geq2$. We have at our disposal ${\mathbb 
V}_{kn}$, $k=1,\dots, i-1$. Here we have that
\begin{equation} \left({}_{Z_{in}}M\right), \ H^o_{nn}\cdots
H^o_{i-1,n}(H^o_{i,n})^{\gamma_X}D^o_{i,n} \in {\mathcal D}(\partial_X,{}_ZM), 
\end{equation} where (use
Lemma~\ref{tedi} in its equivalent $\left({}_{Z_{in}}\right)$ guise)
\begin{equation}\left({}_{Z_{i,n}}M\right)=O=H^o_{i1}\dots
H^o_{i,n-1}H^o_{1n}\dots H^o_{i-1,n}M^o_{in}+\tilde
O=V_{in}^{NW}M^o_{in}+\tilde O.\end{equation} Observing that 
the operators ${\mathbb V}_{kn}$ are monomials in the 
operators $V_{\ell,n},\ell\leq k$, we notice that we could
use the operators $ V_{kn}, k=1,\dots, i-1$ as in
Proposition~\ref{goodprop} to remove the elements in
$\tilde O$ such that the element $V_{in}^{NW}M^o_{in}$
is obtained. Specifically, one could easily use
$V_{1,n},V_{2,n}, \dots $ to eliminate terms $O_j$ for
which $V_{k,n}O_j\neq q O_j V_{k,n}$ (here we must
insist that $\alpha_J\neq2$ as in (\ref{bias}). This
procedure works because each relevant row has a pair
$D,M$ as remarked after Lemma~\ref{tedi}. First use
$V_{1,n}$ then use $V_{2,n}$ on the remaining elements
not distinguished by $V_{1,n}$, then $V_{3,n}$ and so forth. The difference
between the elements $V_{a,n}$ and the (correct)
elements ${\mathbb V}_{an}$ is thus inconsequential,
hence the latter elements also work. Hence $\left(\prod_{s<i}{\mathbb
V}_{sn}\right)V^{NW}_{i,n}M^o_{in} \in {\mathcal
D}_X(\partial_R, LM)$. We also  have the element
$\widehat{RM^-_{W_{1,n}}}^X=H^o_{nn}\dots
H^o_{2,n}(H^o_{in})^{\gamma_X}D^o_{1,n}$, and using
Lemma~\ref{lem-pm} the induction step is completed. \qed

 \medskip

  In complete analogy we get

\begin{Lem}${\mathbb V}_{1,j}$ belongs to
${\mathcal D}(\partial_X,{}_ZM)$ for $j=n,n-1,\dots, 1$. Indeed,
$\left(\prod_{b=j+1}^{n}{\mathbb
V}_{i,b}\right)V^{NW}_{1,j}M^o_{j,1} \in {\mathcal D}(\partial_X,{}_ZM)$ and 
$\widehat{M_{W_{1,j}}}^X\in
{\mathcal D}(\partial_X,{}_ZM)$ where
$M_{W_{1,j}}=V_{1,j}^{SE}M^o_{1,j}$. \end{Lem}

 \medskip

\begin{Prop} For all $i,j$, 
$V_{ij}^{NW}\left(\prod_{a=1}^{i-1}{\mathbb
V}_{aj}\right)M^o_{ij}\in {\mathcal D}(\partial_X,{}_ZM)$ and
$V_{ij}^{SE}\left(\prod_{b=j+1}^{n}{\mathbb
V}_{aj}\right)(H^o_{ij})^{\gamma_X}D^o_{ij}\in {\mathcal D}(\partial_X,{}_XM)$. 
As a consequence, ${\mathbb V}_{ij}\in
{\mathcal D}(\partial_X,{}_ZM)$. \end{Prop}

\proof This follows analogously while using the already established results. We 
can use
$V_{1,j},V_{2,j}\dots, V_{i-1,j}$ for $\left({}_{Z_{ij}}M\right)$ and
$V_{i,j+1},\dots,V_{i,n-1},V_{i,n}$ for $M_{W_{ij}}$.
\qed
 
Remark: We could also use
$V_{i1},V_{i,2},\dots,V_{i,j-1}$ for $\left({}_{Z_{ij}}M\right)$ and, 
independently, 
$V_{i+1,j},\dots,V_{n-1,j},V_{n,j}$ for $M_{W_{ij}}$.
 
 \bigskip

\subsubsection{Right - Right}$\quad$\newline Here, in the terminology of 
Proposition~\ref{goodprop}, ${\mathcal A}_0={\mathcal
D}(\partial_X, M_Z)$, and we have the operators $M_{Z_{ij}}$,  and
$\widehat{M_{W_{ij}}}^X.$ We obtain a result analogous to the previous case, 
but for reasons that should become clear, we find it more reasonable to treat 
the 3 cases one by
one. 

The case J is quite complicated and is considered later. However, we get a 
clean result for each of the cases K, L:
 
\begin{Prop}We have the following: \begin{eqnarray}\nonumber
{\mathcal D}(\partial_D, M_Z)&=&{\mathcal W}eyl_q(n,n) \\\nonumber
{\mathcal D}(\partial_K, M_Z)&=&{\mathbb
C}[V_{ij}^{SE}{\mathcal D}_{ij}\mid 1\leq i,j\leq
n]\cup {\mathbb C}[V_{ij}^{SE}M^o_{ij}\mid 1\leq i,j\leq
n].\end{eqnarray} In particular, $(H^o_{ij})^{\pm2}\in
{\mathcal D}(\partial_K, M_Z)$ for $(i,j)\neq(1,1)$, but only
$(H^o_{11})^{-2}\in {\mathcal D}(\partial_K, M_Z)$.\label{7.4}
\end{Prop}

\pof The induction again follows the ordering \begin{equation}
(n,n)<(n,n-1)< \dots<(n,1)<(n-1,n)<\dots<(1,1).\label{ord}
\end{equation} and starts by observing that $M^o_{nn}$
and $(H^o_{nn})^{\gamma_X}D^o_{nn}\in {\mathcal
D}(\partial_X, M_Z)$. This immediately gives $(H^o_{nn})^{\pm1}$
in case H, while it gives $(H^o_{nn})^{-2}$ and $(H^o_{nn})^{0}$
in case K. To finish case K, we observe that likewise,
$H^o_{n,n}M^o_{n,n-1}$ and $H^o_{nn}(H^o_{n,n-1})^{\gamma_X}
D^o_{n,n-1}\in {\mathcal D}(\partial_X ,M_Z)$. This gives that
$(H^o_{nn})^2(H^o_{n,n-1})^{-2}\in {\mathcal
D}(\partial_K, M_Z)$ and $(H^o_{nn})^2\in 
{\mathcal
D}(\partial_K, M_Z)$. After that the induction proceeds easily
along the same lines with the exception of the restriction in case K
when we reach the final point $(1,1)$. Case L
is the easiest since if we have all $(H^o_{ij})^{\pm1}$
below $(i,j)$ then easily, $M^o_{ij},D^o_{ij}\in
{\mathcal D}(\partial_D, M_Z)$ hence the result. \qed

\medskip

Case $X=J$ is as complicated as in the previous Subsubsection~\ref{rl}, and the 
result is more complex:

\begin{Prop} At the general position $i,j$,
$$V_{ij}^{SE}\left(\prod_{a=i+1}^{n}{\mathbb
Y}_{aj}\right)M^o_{ij}\textrm{ and }
V_{ij}^{SE}\left(\prod_{a=i+1}^{n}{\mathbb
Y}_{aj}\right)H_{ij}D^o_{ij}$$belong to ${\mathcal D}J(\partial_J, M_Z)$; ${\mathbb Y}_{nj}=\prod_{a={j+1}}^n(H^o_{an})^2$ and
${\mathbb Y}_{ij}=(V_{ij}^{SE})^2
\left(\prod_{a=i+1}^{n}{\mathbb Y}_{aj}\right)^2$. 
In particular, ${\mathbb Y}_{ij}\in$ ${\mathcal
D}_J(\partial_J, M_J)$. \end{Prop}

\pof This is again by induction using the ordering (\ref{ord}).

To begin with, at $(n,n)$, we have $M^o_{nn}$ and $H^o_{nn}D^o$.
This leads to $I$ and $(H^o_{nn})^2$. Then we get
$H^o_{nn}M^o_{n,n-1}$ and $H^o_{nn}H_{n,n-1}D^o_{n,n-1}$.
This leads to $(H^o_{nn})^2(H^o_{n,n-1})^2$ (and, once again,
$(H^o_{nn})^2$.) After that it is clear that we have
$(H^o_{nn})^2 (H^o_{n,n-1})^2\dots (H^o_{nr})^2$ for any $n>r\geq1$.
At a position $n-1,j$ we have
\begin{equation}M_{Z_{n-1,r}}=V_{n-1,r}^{SE}M_{n-1,r}^o
+ LOTs\end{equation}and we can use $H_{nn}^2
(H^o_{n,n-1})^2\dots (H^o_{n,r+1})^2$ to separate off the LOTs
terms. A similar result is obtained for
$\widehat{M_{W_{n-1,r}}}^J$ and we then obtain
elements of the form \begin{equation} ((H^o_{nn})^4
(H^o_{n,n-1})^4\dots (H^o_{n,r+1})^4)(H^o_{n,r})^2((H^o_{n-1,n})^2
(H^o_{n-1,n-1})^2\dots (H^o_{n-1,r+1})^2). \end{equation}
 
This element together with the previous can now be used to 
attack a
position $n-2,r$ using Observation~\ref{obs2} and Proposition~\ref{goodprop}. We
use elements with 0 in the upper left corner since they do not affect the top 
term, and in general, we need  one element per row. The result follows. \qed
 
\medskip

 \begin{Rem}This inductive formula above is easily solved: \begin{eqnarray} 
\quad{\mathbb Y}_{ij}&
 =&V_{ij}^{SE}\left(\prod_{x=1}^{n-i}
(Y_{i+x,j}^{SE})^{4\cdot 3^{x-1}}\right)\Rightarrow\\\quad
{\mathbb Y}_{ij}&=&\prod_{x=2}^{n-i}(H^o_{i+x,j})^{2\cdot 
3^{x-2}}(H^o_{i+1,j})^2\prod_{y=j+1}^n (H^o_{iy})^2 \prod_{y=j+1,
x=1}^{n, n-i}(H^o_{i+x,y})^{4\cdot 3^{x-1}}. \end{eqnarray}
\end{Rem} 

We further remark that one can equally well use constructions based on
columns instead of rows.

 \bigskip

\bigskip

\subsection{The algebras of differential operators. Conclusion}

We have seen that the three cases $X=J,L,K$ lead to quite
different results for ${\mathcal D}(\partial_X, M_X)$ and
${\mathcal D}(\partial_X, {}_ZM)$ though the following  holds that in all 6 
($=3\times 2$)
cases:
\begin{Thm}For any index $i,j$ there are, case dependent,
elements $\psi_{ij}$ and $\phi_{ij}$ of ${\mathcal
L}{\mathcal M}$ such that \begin{equation}\psi_{ij}M^o\textrm{  and 
}\psi_{ij}D^o\end{equation} belong to the algebra of the given case.
\end{Thm}

\begin{Def}Set\begin{equation}{\mathcal D}^R_X ={\mathcal
D}(\partial_X, M_Z,{}_ZM); \end{equation}the algebra of quantized right 
differential
operators, and set\begin{equation}{\mathcal D}^L_X={\mathcal
D}({}_X\partial,M_Z,{}_ZM); \end{equation}the algebra of quantized left 
differential
operators.\end{Def}

We now discuss the cases J, K, L one by one:

\begin{Thm}\label{L-thm} We have that \begin{equation}{\mathcal D}^R_L 
={\mathcal
D}^L_L ={\mathcal W}eyl_q(n,n). \end{equation} The representation ${\mathcal 
O}_{\Lambda'}^L$ is given by differential operators together with 
 $\End(V'_\Lambda)$.
\end{Thm}

\pof It is clear that \begin{equation}{\mathcal
D}({}_H\partial,{}_ZM)={\mathcal D}(\partial_H, M_Z)={\mathcal
W}eyl(n,n). \end{equation}The rest is equally obvious. \qed

\medskip

In case J we do not get $\widehat{\left({}_{W_{11}}M\right)}^J$ from the
right and we do not get $\widehat{M_{W_{nn}}}^J$ from the
right. We do get parts of $\widehat{\left({}_{Z_{11}}M\right)}^J$ (on $Z$s) 
but miss out on  eg. $K_\beta^{-1}$.
 We will refrain from stating further results for this case.

\medskip

\begin{Def} \ ${\mathcal K}{\mathcal W}eyl_q(n,n)$ is set to be the smallest 
sub\-alge\-bra
of  ${\mathcal W}eyl_q(n,n)$ containing the following operators: 
\begin{equation}
\forall i,j: V_{ij}^{NW}(H^o_{ij})^{-1}D^o_{i,j},
V_{ij}^{NW}M^o_{i,j},\textrm{ and }
V_{ij}^{SE}M^o_{i,j}. \end{equation} In particular, it
contains all $(H^o_{ij})^{\pm2}$ and
$V_{ij}^{NW}V_{ij}^{SE}$.
\end{Def}
\medskip

\begin{Thm}\label{K-thm} We have that \begin{equation}{\mathcal D}^R_K 
={\mathcal
D}^L_K:= {\mathcal K}{\mathcal W}eyl_q(n,n).\end{equation} The representation 
${\mathcal O}_{\Lambda'}^K$ is given by operators from ${\mathcal K}{\mathcal 
W}eyl_q(n,n)$ together with $End(V'_\Lambda)$.
\end{Thm}

\pof The fact that the left and right algebras are
identical is easy to see since by symmetry we easily get  ${\mathcal
D}({}_K\partial_L,{}_ZM)$ while ${\mathcal
D}({}_K\partial,M_Z)$  is the dual of an analogous algebra
 ${\mathcal
D}_K({}_WM,{}_K\partial )$ on the ${\mathcal A}_q^-$. The statements
about the elements in the algebra follow from
Proposition~\ref{7.4} and the following observation: The combined
operators easily lead to e.g. $H^o_{11}V^{SE}_{12}\in
{\mathcal D}^R_K$. By squaring this element we easily
obtain $(H^o_{11})^2$.

As for the the representation, by definition
$\widehat{L_{W_\beta}}^K\in {\mathcal D}_K^R$. Next
observe that $K_\beta=(H^o_{11})^2V^{11}{NW}V^{SE}\in
{\mathcal D}_K^R$ and similarly, since we have all
operators $(H^o_{ij})^{-2}$ and thus can replace odd
positive powers in elements $H^o_{st}$ by odd negative
powers, $K^{-1}_\beta\in {\mathcal D}_K^R$.  Now consider ${\mathcal 
O}^K_{\Lambda'}(E_\beta)$: It differs from ${\mathcal 
O}^J_{\Lambda'}(E_\beta)$ by squares of elements $H_{ij}^o$, cf. 
(\ref{h1}-\ref{h3}). We may therefore equally well look at the latter which, by 
Corollary~\ref{actJ} is given by by $K^{-1}_\beta$  and left and right 
multiplication operators. We remark here that for the general expression in 
Corollary~\ref{actJ} we need to introduce $End(V'_\Lambda)$ (which is generated 
by the operators from ${\mathcal U}_1({\mathfrak k}^{\mathbb C})$. We then only 
need to consider the representation restricted to ${\mathcal U}_1({\mathfrak 
k}^{\mathbb C})$, and we can here, by similar arguments, restrict to consider 
the scalar case and we may as well just consider the expressions in 
Proposition~\ref{k-prod} and easily get, up to factors of $(H^o_{st})^{\pm2}$,

\begin{eqnarray} E_{\mu_k}^Z&=&
\sum_j(V_{k,j}^{NW}(H^o_{k,j})^{-1}
D^o_{k,j}V_{k+1,j}^{NW}M^o_{k+1,j}),\textrm{ and}\\ F_{\mu_k}^Z&=&
\sum_j(V_{k+1,j}^{SE}(H^o_{k+1,j})^{-1}
D^o_{k+1,j}V_{k,j}^{SE}M^o_{k,j}). \end{eqnarray}

Since these expressions are in the algebra, we are done. \qed

\bigskip

We finish with the following result which of course is immediate from the way 
these alebras are invariantly defined from left and right actions:

\begin{Thm}  ${\mathcal W}eyl_q(n,n)$ as well as the various constructed 
subalgebras, though given  manifestly in a fixed PBW basis, are intrinsically 
invariant.
\end{Thm}


\bigskip



\begin{thebibliography}{22}

\bibitem{chari} V. Chari, A. Pressley. A Guide to Quantum Groups. Cambridge 
Univ. Press, 1995.

\bibitem{hi} T. Hayashi, $q$-analogues of Clifford and Weyl Algebras--Spinor 
and Oscillator Representa\-tions of Quantum Enveloping Algebras, Comm. Math. Phys. Volume 127, Number 1 
(1990),
129-144.

\bibitem{jv} Hans P. Jakobsen and Michele Vergne, Restrictions and expansions 
of holomorphic representations, J. Functional Analysis 34, 29-53 (1979)

\bibitem{j-qdiff} Hans P. Jakobsen, Q-differential operators. Preprint (19pp) 
1999\newline (http://xxx.lanl.gov/abs/math.QA/ 9907009)

\bibitem{j-qdirac} Hans P. Jakobsen, Quantized Dirac operators. Czech. J. Phys. 
50, 1265--1270 (2000)



\bibitem{jz} H. P. Jakobsen and H. Zhang, Double-partition Quantum Cluster
Algebras J.
Algebra 372 (2012), 172-203. 


\bibitem{jan} J.C. Jantzen, {\em Lectures on Quantum Groups},
A.M.S. Graduate 
Studies in Mathematics Vol. 6 (1996)

\bibitem{jo-le} A. Joseph and G. Letzter, Rosso's Form and Quantized Kac Moody 
algebras, Mathematische Zeitschrift, Vol. 222 (1996) 543-571.


\bibitem{K1} M. Kashiwara, On crystal bases of the $q$-analogue of universal
enveloping algebras,
Duke Math. J. 63 (1991), 465--516.

\bibitem{K2} M. Kashiwara, Global crystal bases of quantum groups, 
Duke Math. J.  69 (1993), 455--485.

\bibitem{kilu}  A. P. Kitchin and S. Launois, On the Automorphisms of Quantum 
Weyl Algebras, arXiv:1511.01775v2 [math.QA] 12 Jul 2018, to appear in Journal 
of Pure and Applied Algebra.


\bibitem{lu1} G. Lusztig, Quantum deformations of certain simple modules
over enveloping algebras,  Adv. in Math. 70 (1988),  237-249. 

\bibitem{luz}G. Lusztig, Canonical bases arising from quantized enveloping 
algebras, J Amer. Math. Soc. 3 (1990), 447-498.

\bibitem{luz2}G. Lusztig, Quantum groups at a root of 1,
 Geom. Dedicata, 35, 1990, p. 89-144.

\bibitem{l} G. Lusztig {\bf Introduction to Quantum Groups.}
Progress In
Mathematics 110, Birkh{\"a}user 1993.

\bibitem{rosso}M. Rosso, Certaines formes bilinéaires sur les groupes 
quantiques et
une conjecture de Schechtman et Varchenko, 
C. R. Acad. Sci. Paris, t. 314, S\'erie I, p. 5-8, 1992

\bibitem{vaks} S. Sinel'shchikov, L. Vaksman, Hidden symmetry of the 
differential calculus on the quantum
matrix space, J. Phys. A., Math. Gen., 30 (1997), 23--26; in Lectures on 
q-analogs of Cartan
domains and associated Harish-Chandra modules, L. Vaksman (ed.), 
math.QA/0109198,
136--140.

\bibitem{genkai} D. Shklyarov, G. Zhang, Covariant $q$-differential operators and unitary highest weight 
representations for $U_q(su(n,n))$, J. Math. Phys. {\bf46} Issue 6 (2005). 	
dx.doi.org/10.1063/1.1927077

\bibitem{tani} T. Tanisaki, Killing forms, 
Harish-Chandra Isomorphisms, and universal R-matrices for quantum algebras, pp. 
941-961 in: A. Tsuchiya, T. Eguchi, and M. Jimbo (eds.), {\bf Infinite 
Analysis}, Part B, Proc. Kyoto 1991 (Advanced Series in Mathematical Physics 
16), Rivers Edge, N. J., 1992 (World Scientific)


\end{thebibliography}
\end{document}